\documentclass[preprint]{elsarticle}

\usepackage{lineno,hyperref}
\modulolinenumbers[5]

\usepackage[utf8]{inputenc}
\usepackage{amsmath}
\usepackage{amssymb}
\usepackage{fullpage}
\usepackage{graphicx}
\usepackage{amsthm}
\usepackage{rotating}
\usepackage{float}
\usepackage[intoc]{nomencl}
\makenomenclature
\usepackage{natbib}
\usepackage{setspace}
\usepackage{listings}
\usepackage{color}
\usepackage{gnuplot-lua-tikz}
\usepackage{tikz}
\usetikzlibrary{positioning,shapes,arrows,backgrounds,calc,fit}
\usepackage{parskip}
\usepackage{subfigure}
\usepackage[nottoc]{tocbibind}
\usepackage{footnote}
\makesavenoteenv{minipage}
\makesavenoteenv{itemize}
\usepackage{hyperref} 

\newcommand{\D}[2]{\frac{\partial #1}{\partial #2}}

\renewcommand{\vec}[1]{{\bold #1}}

\newcommand{\bm}[1]{\mathbf{#1}}

\newcommand{\firstpktot}{\hat{\bm{P}}}
\newcommand{\firstpkref}{\bm{P}}
\newcommand{\secondpktot}{\hat{\bm{S}}}
\newcommand{\secondpkref}{\bm{S}}

\newcommand{\heatfluxref}{\vec{q}}
\newcommand{\deftentot}{\bm{F}}

\definecolor{mygreen}{rgb}{0,0.6,0}
\definecolor{mygray}{rgb}{0.5,0.5,0.5}
\definecolor{mymauve}{rgb}{0.58,0,0.82}

\theoremstyle{definition}

\theoremstyle{remark}

\journal{Finite Elements in Analysis and Design}

\def\statement{\begin{minipage}[t]{.75\textwidth}
       NOTICE: This is the author's version of a work that was accepted for
publication in Finite Elements in Analysis and Design.  Changes
resulting from the publishing process, such as editing, corrections, structural
formatting, and other quality control mechanisms may not be reflected in this
document. Changes may have been made to this work since it was submitted for
publication.
\\

\copyright \, 2019. This manuscript version is made available under the
CC-BY-NC-ND 4.0 license \\
\url{http://creativecommons.org/licenses/by-nc-nd/4.0/}.
       \end{minipage}}

\makeatletter
\def\ps@pprintTitle{%
     \let\@oddhead\@empty
     \let\@evenhead\@empty
     \def\@oddfoot{\footnotesize\itshape
       \statement\hfill\today}%
     \let\@evenfoot\@oddfoot}
\makeatother

\journal{}

\begin{document}

\begin{frontmatter}

\title{Inverse Design Based on Nonlinear Thermoelastic Material Models Applied
to Injection Molding}

\author{Florian Zwicke\corref{mycorrespondingauthor}}
\cortext[mycorrespondingauthor]{Corresponding author}
\ead{zwicke@cats.rwth-aachen.de}

\author{Stefanie Elgeti}
\ead{elgeti@cats.rwth-aachen.de}

\address{Chair for Computational Analysis of Technical Systems,
RWTH Aachen University, Schinkelstr. 2, 52062 Aachen, Germany}

\begin{abstract}
This paper describes an inverse shape design method for thermoelastic bodies.
With a known equilibrium shape as input, the focus of this paper is the
determination of the corresponding initial shape of a body undergoing
thermal expansion or contraction, as well as nonlinear elastic deformations.
A distinguishing feature of the described method lies in its capability to
approximately prescribe an initial heterogeneous temperature distribution
as well as an initial stress field even though the initial shape is unknown.
At the core of the method, there is a system of nonlinear partial differential
equations. They
are discretized and solved with the finite element method or isogeometric
analysis.
In order to better integrate the method with application-oriented simulations,
an iterative procedure is described that allows fine-tuning of the results.
The method was motivated by an inverse cavity design problem in injection
molding applications. Its use in this field is specifically highlighted,
but the general description is kept independent of the application to simplify
its adaptation to a wider range of use cases.
\end{abstract}


\end{frontmatter}

\section{Introduction}

One of the principal strengths of numerical simulations is that they can give
insights into the behavior of products before they are actually built. These
insights can even be useful in improving product design before the first
prototype is produced.
In a best-case scenario, numerical simulations could actually be used to find
the perfect product design on the computer, without having to produce any
prototypes.

In this paper, we deal with inverse design problems for thermoelastic materials
with nonlinear elastic behavior.
Our goal is to determine an initial shape of a body, when all other relevant
parameters, as well as the body's shape under equilibrium conditions, are known.
These types of problems are also called shape optimization problems.


In order to solve these types of problems, we have developed a method
that is general enough to be
applicable to a broad range of applications, mostly in the field of production
engineering.
We will specifically
highlight the method's usefulness in the field of injection molding, for which it has
been developed. However, in the way the method is presented, or with just minor
adjustments, it can also be applied to processes such as
high pressure die casting or additive manufacturing (e.g., selective laser
melting).

\begin{figure}[htbp]
  \centering
  \includegraphics[height=2.5cm]{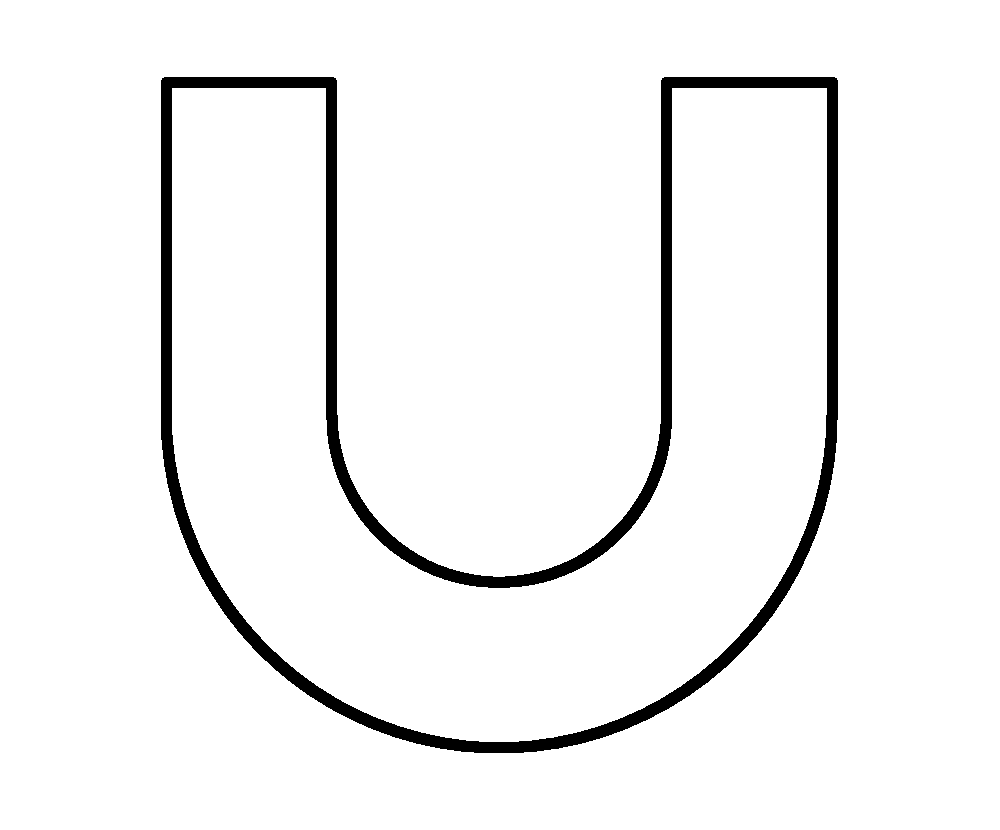}
  \begin{minipage}[b]{1cm}
    {\Huge $\,\rightarrow\,$}

    \vspace{0.8cm}
  \end{minipage}
  \includegraphics[height=2.5cm]{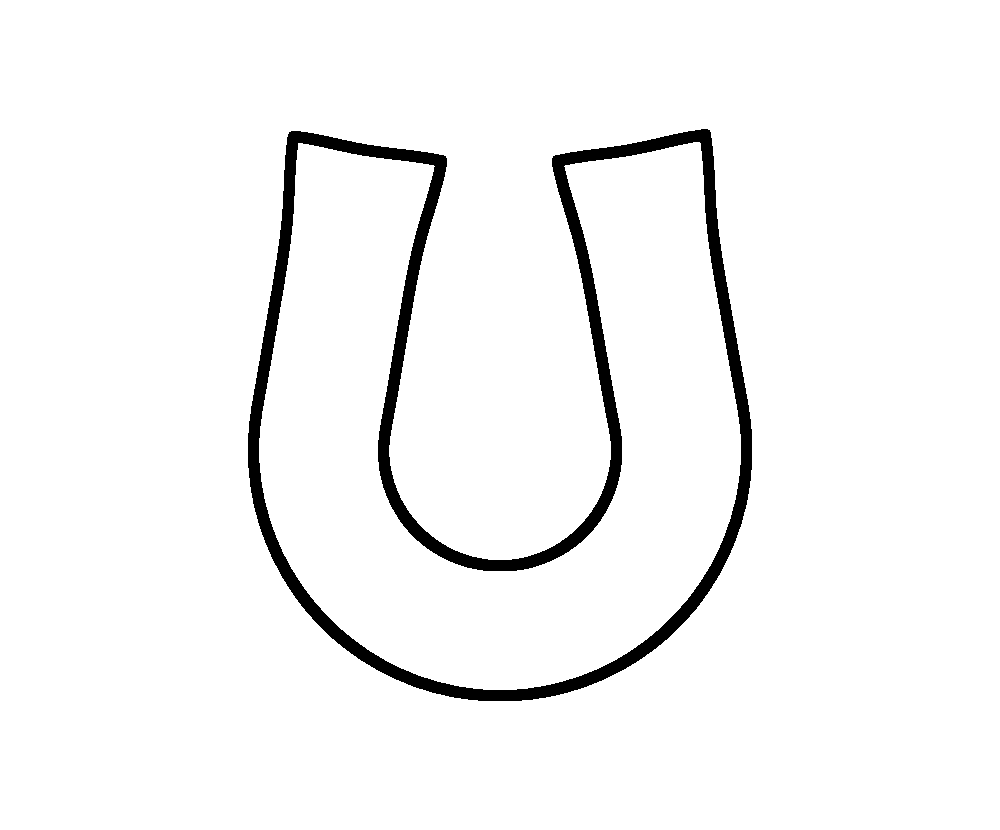}
  \caption{When a certain cavity shape is used (left-hand side), a certain
molding shape results (right-hand side).}
  \label{fig:introforward}
\end{figure}

Injection molding is a process where a liquid polymer melt is injected into a
cavity where it is cooled down such that it solidifies \cite{Potsch2008}.
Both the cooling and solidification cause a decrease in the specific volume of
the material. Since this volume change happens inhomogeneously, the shape of the
body changes such that the end result is not just a scaled version of the cavity shape.
Figure~\ref{fig:introforward} is a simple sketch of how a molding's shape could
deviate from the cavity shape.

\begin{figure}[htbp]
  \centering
  \includegraphics[height=2.5cm]{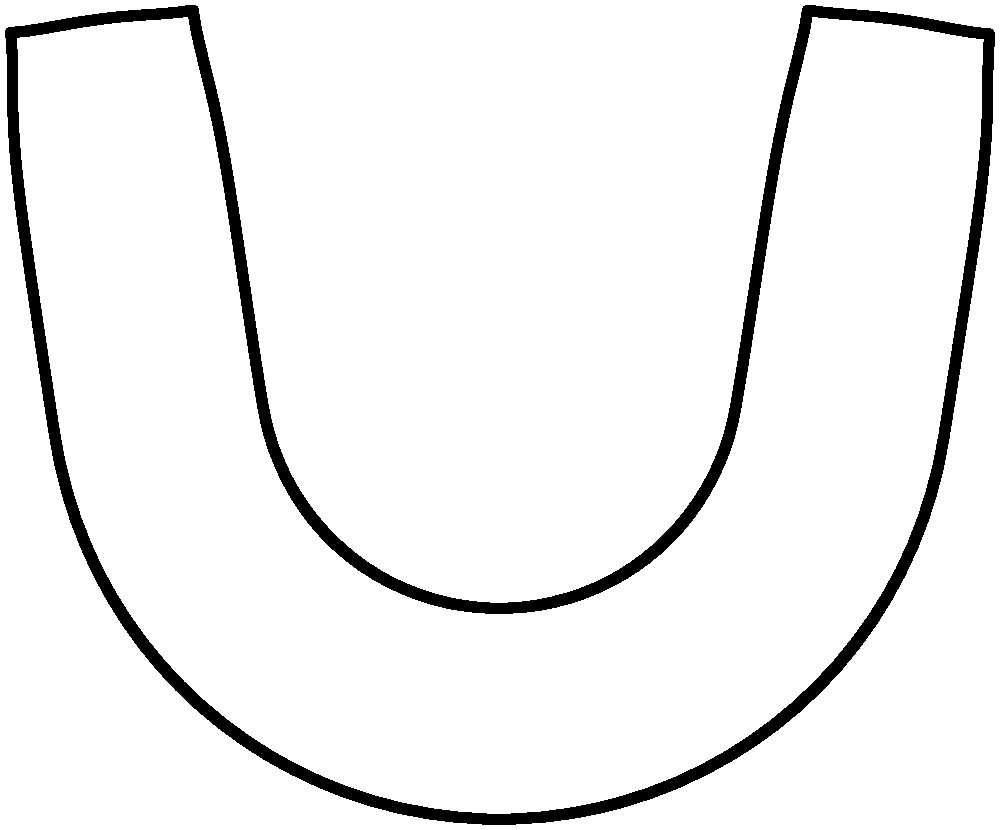}
  \begin{minipage}[b]{1cm}
    {\Huge $\,\leftarrow\,$}

    \vspace{0.8cm}
  \end{minipage}
  \includegraphics[height=2.5cm]{fig-shape-desired.png}
  \caption{The inverse method we seek can generate a cavity shape (left-hand side), when
a molding shape is prescribed (right-hand side).}
  \label{fig:introreverse}
\end{figure}

Of course, we would like to be able to prescribe the molding's shape exactly.
Thus, what we are looking for is an inverse method that can yield a cavity
shape that produces the correct molding shape.
This is illustrated in Figure~\ref{fig:introreverse}.

The injection molding process can be subdivided into several stages. In this
work
we will focus on the part of the process that happens after the material has
solidified enough to be considered fully elastic.
This means that our simulation starts from a point in time when there is a solid body that
has a shape that is still identical to the cavity shape. This body has
inhomogeneous temperature and stress distributions
-- these can for example be determined from a filling and solidification
simulation or also from experiments.
The body eventually cools down to the environment temperature, and the
stresses relax to a point where the body is in equilibrium, i.e., the
normal stresses on the boundary become zero. During this process, the body
undergoes significant shape changes.
It is the simulation of this process that we have to invert to determine a useful
cavity shape from a prescribed molding shape.
For a detailed description of how this method can be applied to injection
molding and what other approaches are available, we refer to~\cite{Zwicke2017}.
To treat the problem in a more
general fashion, independent of the specific application of injection molding,
we will now refer to the cavity shape as the body's \textit{initial shape}, and
the molding shape as the \textit{equilibrium shape}.


In general, when a forward simulation is available, i.e., a simulation that can predict the
body's equilibrium shape given an initial shape, this problem can be dealt with
in the context of mathematical optimization. One further requirement for this is
the availability of an objective function, i.e., a scalar measure that can describe
the geometric differences
between some output shape of the simulation and the desired equilibrium shape.
This would make it possible to use, e.g., methods of quasi-Newton type for the
optimization, such as BFGS \cite{Liu1989}, or derivative-free optimization methods.
Optimization methods have been used in the past to
solve problems very similar to ours \cite{Nobrega2008, Ettinger2004, Elgeti2012, Siegbert2015}.


Optimization methods can be used as soon as a working forward
simulation for the problem is available. This makes them very versatile, but
they also present challenges.
Most importantly, quasi-Newton methods can only be applied if the derivative of
the objective function with respect to the design parameters is available. In
shape optimization problems, these parameters describe the unknown shape.
Therefore, their number can be quite large for complex shapes.
In cases, such as ours, where the objective function depends on the result of a
rather complex simulation, the calculation of the derivative is difficult. This
means that it either has to be approximated, or derivative-free
optimization methods have to be used. Both options are time-consuming for
large numbers of design parameters.

The optimization problem we are dealing with is a very special type of shape
optimization problem, where the objective function is also a shape.
Additionally, in contrast to other optimization problems, where the goal is to
minimize some objective function without prior knowledge of the simulation
result or minimum value, we actually know the result that we wish to obtain.
This means that what we are trying to achieve is an actual inversion of the
simulation.
Such an inverse simulation can be much more efficient than the aforementioned
shape optimization methods. Furthermore, this can remove the need to
parameterize the unknown shape --- a step that limits the design space to
less-than-optimal solutions.


For isothermal elasticity, inverse formulations have been researched already for
some time.
In 1967, Schield~\cite{Schield1967} showed that an inverse
nonlinear elasticity problem, where the equations are solved for the reference rather than
the deformed state, has the same general form as a regular (i.e., forward)
elasticity problem.
This was later extended by Carlson~and~Schield~\cite{Carlson1969}.
Govindjee~et~al.~\cite{Govindjee1996, Govindjee1998} showed how the standard equilibrium equations can
be re-parameterized to formulate the inverse problem. In this case, the
Lagrangian viewpoint typical of elastostatics is exchanged by a Eulerian
viewpoint, where the known equilibrium configuration is used as the reference
and the unknown initial configuration is solved for.
\cite{Hong2016} is another example of a similar formulation.

As we will point out in Section~\ref{sec:solstrat}, we require the reference
configuration to be different from both the equilibrium and initial
configurations. Such an ALE-type (Arbitrary Lagrangian-Eulerian) approach
was used by Yamada in 1998, where this was
applied to incompressible hyperelasticity \cite{Yamada1998}.

Such inverse formulations have been used, among other things, for the
optimization of turbine blade shapes \cite{Bazilevs2012, Campbell2011, Fachinotti2008},
large beams \cite{Albanesi2008}, compliant mechanisms \cite{Albanesi2013, Albanesi2011,
Albanesi2009}, and
airplane wings \cite{Limache2011}.


In contrast to many inverse elasticity problems, the main source of deformations
are, in our case, not external forces but temperature changes, i.e., thermal
expansion or contraction. The notion of inverse formulations in thermoelasticity
has appeared in literature before, although sometimes in a slightly different
context.
Dennis~et~al. have applied the term to problems where the initial state is
known, but parts of the boundary
conditions in the equilibrium state are unknown,
i.e., neither displacements nor forces are known on parts of the boundary.
This is made up for by over-specifying boundary conditions on other parts of the
boundary, which requires both displacements and forces to be prescribed
\cite{Dennis2004, Dennis2011}.
This means that both Dirichlet and Neumann boundary conditions are applied
simultaneously in the same place in some cases.
Such an over-specification of Dirichlet and Neumann boundary conditions also happens in our case,
but on the entire boundary, to allow the prescription of an equilibrium shape
at the same time as zero normal stresses on the boundary.
The major difference to the work of Dennis~et~al. lies in the fact that we use
this over-specification to solve for displacements in the initial state, which
creates the need for multiple displacement fields.
This is explained in
more detail in Section~\ref{sec:invsim}.

\section{Solution Strategy}
\label{sec:solstrat}

\subsection{Description of the Inverse Problem}

Our problem has a slightly different setup from many other inverse
design problems in elastostatics. In many cases, two states of a body
are considered. There is one state where the body is subjected to external
loads and body forces such as gravity. In a second state, the body is
considered without any forces acting on it.
We have a distinct problem where a body
transitions from a constrained to an unconstrained state. We will call
the chronologically first state the \textit{initial state}, and the latter
state the \textit{equilibrium state}, since the body will be in thermodynamic
equilibrium.

Before we detail the setup of the inverse problem, we will first fully describe
the corresponding forward simulation. In the initial state, we are dealing with
an inhomogeneously heated body that is physically constrained. The state of the
body is fully described by its shape, the field of internal stresses and the
temperature field. In the equilibrium state, we have knowledge of the
environment temperature, gravitational forces and external forces, which may be
zero in many cases. The principal quantity of interest is the shape of the body
in the equilibrium state. However, the field of internal stresses in this state is also
needed, since the stresses need to stay below a certain threshold for the
elastic material models to remain valid.
Table~\ref{tab:forwardsim} shows an overview of the
prescribed and unknown quantities in both states.

\begin{table}[htbp]
  \centering
  \caption{Forward simulation}
  \begin{tabular}{| p{5cm} || p{5cm} | p{5cm} |}
    \hline
    & \textbf{initial state} & \textbf{equilibrium state} \\
    \hline
    \hline
    \textbf{prescribed quantities} & \begin{itemize}
      \item boundary shape
      \item internal stresses
      \item temperature
      \end{itemize} &
      \begin{itemize} \item temperature \item external and internal forces
        \end{itemize} \\
    \hline
    \textbf{unknown quantities} & & \begin{itemize} \item boundary shape
      \item internal stresses \end{itemize} \\
    \hline
  \end{tabular}
  \label{tab:forwardsim}
\end{table}

The inverse problem results from the wish to prescribe the equilibrium shape of
the body.
One particularity of this inverse problem is the necessity to prescribe a
temperature and stress distribution on an initial shape that is still to be computed.
Table~\ref{tab:inversesim} shows the changed distribution of
prescribed and unknown quantities for the inverse simulation.

\begin{table}[htbp]
  \centering
  \caption{Inverse simulation}
  \begin{tabular}{| p{5cm} || p{5cm} | p{5cm} |}
    \hline
    & \textbf{initial state} & \textbf{equilibrium state} \\
    \hline
    \hline
    \textbf{prescribed quantities} &
      \begin{itemize} \item \textit{internal stresses} \item \textit{temperature} \end{itemize} &
      \begin{itemize} \item boundary shape \item temperature
        \item external and internal forces \end{itemize} \\
    \hline
    \textbf{unknown quantities} & \begin{itemize} \item boundary shape
      \end{itemize} &
      \begin{itemize} \item internal stresses \end{itemize} \\
    \hline
  \end{tabular}
  \label{tab:inversesim}
\end{table}

\subsection{Handling of Auxiliary Fields}
\label{sec:auxfields}

In Table~\ref{tab:inversesim}, the internal stresses and temperature are
italicized since they require special treatment. It is difficult to prescribe,
for instance, an inhomogeneous temperature field, if the shape of the body is
unknown. A mechanism needs to be devised that can produce the correct
temperature field for a certain shape.
An important requirement of this mechanism is that it has to be integrated
seamlessly with the system of equations for the inverse problem to ensure that
it can still be solved in a \textbf{monolithic}, i.e., non-iterative, fashion.

To initialize the inverse simulation, the temperature and stress fields are
determined for a reference configuration chosen by the user (cf. also
Section~\ref{sec:refconf}): In injection
molding, one might perform a filling simulation for a cavity resembling a
scaled-up version of the product shape.
These fields, from now on referred to as reference temperature and reference
stress, will then serve as a basis for estimating the corresponding fields in
arbitrary shapes; a step which is of course only necessary if the fields cannot
be computed directly.
For this estimation, we use a mechanism that is identical
for both the temperature and stress field, and could be applied to other
auxiliary fields. Therefore, we will only describe it for the temperature field.

The temperature field in the initial state entirely depends on the application
that is considered. Even small changes in the boundary shape could mean an
entirely different temperature field. For our application of injection molding,
we make the assumption that small changes in the boundary shape will also cause
small changes in the temperature field. Following this logic, we look at the
difference between a certain shape and the reference shape, and try to translate
this shape difference to a difference in the temperature distribution. The only
requirement for this translation is that for identical shapes, the temperature
fields should also be identical.

As long as we do not wish to incorporate any knowledge of how the temperature
field is produced in the first place into the system of equations,
we should always aim to keep the
difference between the estimated and reference temperature fields as small as
possible, as long as we fulfill the constraint of the adjusted boundary shape.

In order to achieve this, we borrow an idea from interface tracking methods
in free-surface flow simulations. In such contexts, the simulation mesh needs to
be updated to fit a changed boundary in a way that avoids damage to the mesh
in terms of deteriorating simulation properties \cite{Zwicke2017a}. A popular method for this
purpose is called EMUM, the Elastic Mesh Update Method \cite{Johnson1994}.
In this method, an elasticity problem is solved with Dirichlet boundary
conditions to achieve, in a sense, minimal movement of interior mesh nodes.
We apply the same method to move temperature nodes such that they fit a new
shape. This means that we assume no actual changes in temperature values, but
only in their distribution. One should note a small imperfection in this method,
which stems from the fact that tangential sliding of boundary nodes, even if
this does not contribute to an actual shape change, will induce unwanted changes
in the temperature distribution. We currently neglect this issue in favor of
efficiency.
It is possible to use a linear elasticity model for this purpose in cases where
the deformation is small enough. In other cases, more
sophisticated equations, such as nonlinear elasticity, should be considered.
In our description, we will refer to the equation for the node movements
as \textit{smoothing equation}.


\subsection{Reference Configuration}
\label{sec:refconf}

We have so far introduced three different states of a body: the initial state,
the equilibrium state, and a reference state.
We will now explain in more detail the connection between the configurations of
the body in these different states.
For a structural finite element simulation, we need to define one reference
configuration. This defines the simulation mesh as well as the coordinate system
used to formulate the differential equations. We will define the set of points
that belong to the body in the reference configuration as $\Omega \subset
\mathbb{R}^3$. This allows
us to describe the other configurations using displacement fields
$\vec{d} \in \Omega$ that store the displacement vectors of all material points.

In a regular, i.e., forward-in-time, structural simulation,
the initial configuration may be
chosen as the reference. This is convenient, since the initial configuration is
known and the displacement field for the equilibrium configuration
immediately translates to
the unknown quantity that is solved for. This approach is referred to as the
\textit{Lagrangian formulation}.

If we consider an inverse structural simulation, where the equilibrium
configuration is fully prescribed, the opposite approach can be applied. This
means that the equilibrium configuration serves as the reference, and a
displacement field for the initial configuration is calculated. This is called the
\textit{Eulerian formulation}.
As we already mentioned, this approach has been used in many cases
to formulate inverse design problems \cite{Schield1967,Govindjee1996}.

For our simulation, we have already mentioned a reference configuration
that is distinct from both the initial and the equilibrium configurations.
This is the configuration for which the auxiliary fields, such as
temperature and internal stresses, that are used to estimate these fields in the
initial configuration, are provided. We introduce the symbols

\begin{itemize}
  \item $\vec{r}$ for the displacement field describing the initial configuration,
        and
  \item $\vec{u}$ for the displacement field
        describing the equilibrium configuration,
\end{itemize}

both defined on the reference configuration $\Omega$.

The aforementioned smoothing equation for the
transformation of the auxiliary fields determines the interior of the displacement
field $\vec{r} \in \Omega \setminus \partial \Omega$ from the boundary
displacement $\vec{r} \in \partial \Omega$. For this reason, the interior values
of the displacement field $\vec{u} \in \Omega \setminus \partial \Omega$
need to be left unspecified. This means that neither $\vec{r} \in \Omega$ nor
$\vec{u} \in \Omega$ are fully prescribed, so a different reference must be
used. This third approach is called \textit{ALE}, the
Arbitrary-Lagrangian-Eulerian formulation. In the context of inverse design
problems, this formulation can be found in, e.g., \cite{Yamada1998}.

\begin{figure}[htbp]
  \centering
    \input{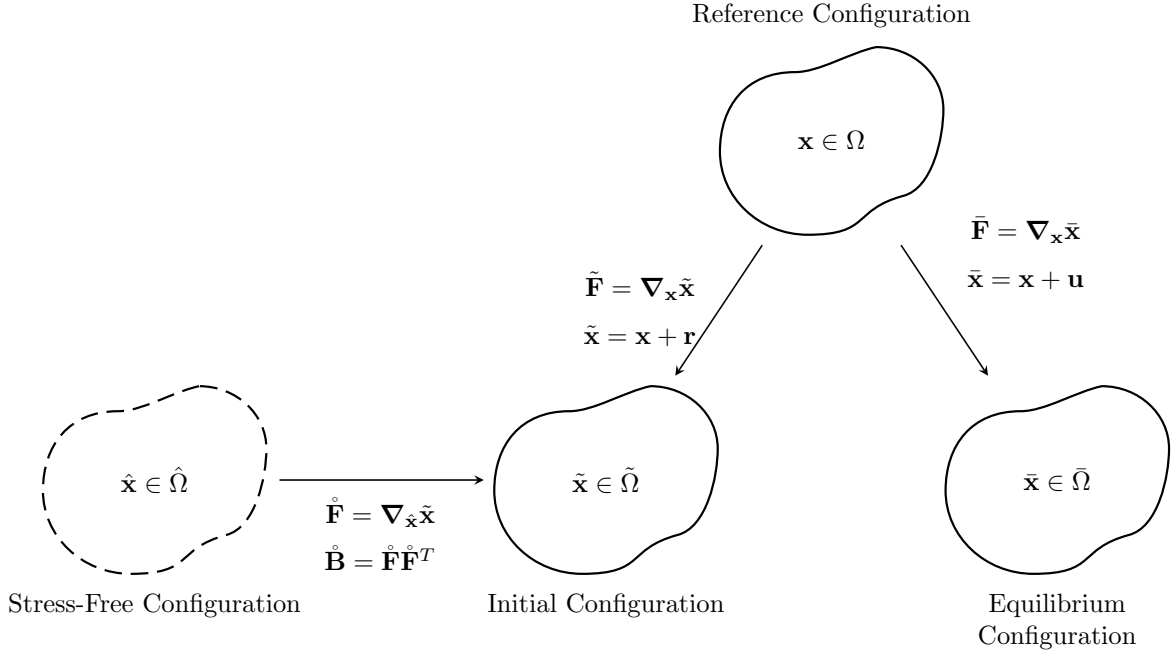}
  \caption{Relations between the different configurations}
  \label{fig:configurations}
\end{figure}

Figure~\ref{fig:configurations} shows a sketch of how the different
configurations relate to each other and which symbols are used to describe these
relations. We have introduced the additional symbol $\mathring{\bm{B}}$ to
describe the deformation from a theoretical, stress-free configuration to the initial
configuration. This field indirectly stores the initial stresses.
It makes sense to keep track of the initial deformation rather than the stress
resulting from it, since this can more easily be integrated with existing
stress-strain relationships. Furthermore, we only require the symmetric tensor $\mathring{\bm{B}}$
instead of $\mathring{\bm{F}}$, since this suffices in connection with the
material laws that we are using and since it allows us to use results from
simulations that do not explicitly keep track of
$\mathring{\bm{F}}$.

\subsection{Iteration Scheme}

The method for handling auxiliary fields that we presented in
Section~\ref{sec:auxfields} only approximates these fields based on an
assumption of minimal change. However, the true effects of the geometry changes
that are proposed by the inverse method on these fields are unknown. If we
intend
to incorporate these effects, we have to rerun the simulations that yield these
auxiliary fields with the changed geometry. Such a rerun will offer new
initial data for the inverse method, so the initial shape can in turn be improved.
Repetition of these steps leads to an iteration scheme, such as the following:

\begin{enumerate}
  \item Select best guess for initial shape (e.g. desired shape).
  \item Determine auxiliary fields for current shape.
  \item Run forward simulation to determine the equilibrium shape.
  \item If the quality of the equilibrium shape is good enough (or another
        stopping criterion, such as maximum number of iterations, is reached), stop the
        iteration.
  \item Run inverse simulation to determine improved guess for initial shape.
  \item Continue at Step 2.
\end{enumerate}

\begin{figure}[htbp]
  \centering
    \input{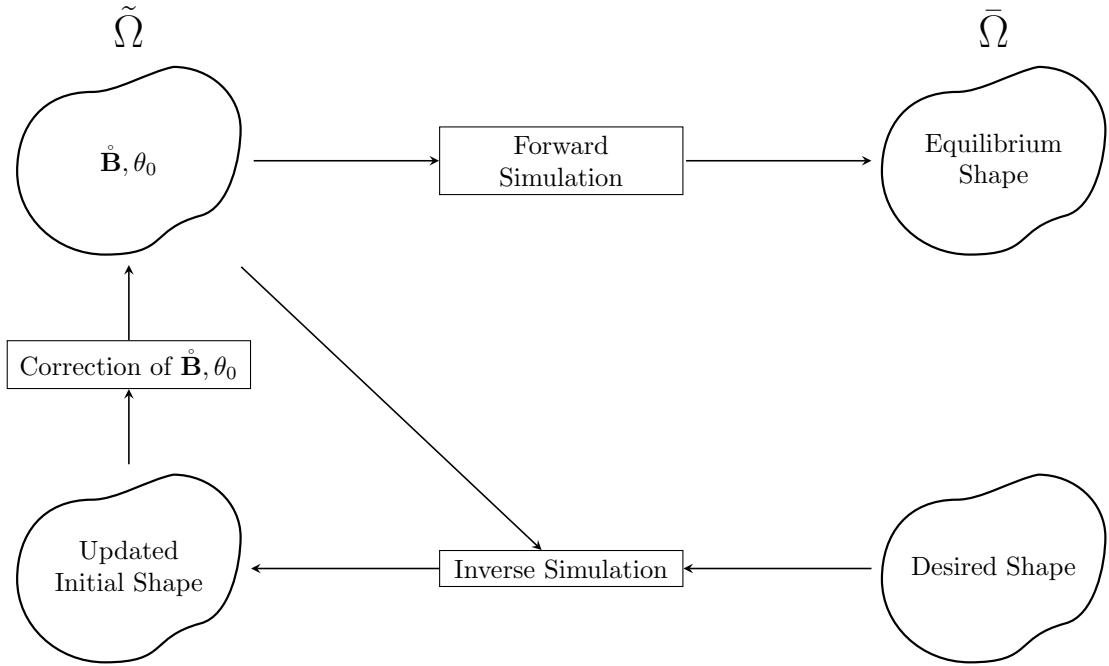}
  \caption{An iteration scheme can combine the inverse simulation and a
correction step for the initial deformations $\mathring{\bm{B}}$ and the initial
temperature $\theta_0$ to counteract the errors introduced by an inexact
smoothing equation.}
  \label{fig:iterationscheme}
\end{figure}

This iteration scheme is illustrated in
Figure~\ref{fig:iterationscheme}. The correction step that is mentioned in the
illustration could be, as previously mentioned, another simulation, or some
other method to determine the auxiliary fields.

There is no guarantee that such an iteration will converge. This completely
depends on how well the true changes in the auxiliary fields match with the
assumptions made in the smoothing equation.

In the specific case of cavity shape design in injection molding, the correction
step consists of a viscoelastic fluid simulation, including polymer
solidification models, that can generate a temperature and stress field for
a provided cavity shape.

\section{Conservation Laws in Regular and Inverse Formulations}

After having described the inverse problem, we will now describe the equations
that are needed to solve it. These include the equilibrium equation for the
elastic body as well as the smoothing equation for the initial displacement
field. We will also show how these equations can be discretized to be solved
using the finite element method or isogeometric analysis.

\subsection{Smoothing Equation}

The equation that is needed to find an approximation for the changes in the
auxiliary fields should lead to minimal changes in the interior of the field to
match the changes on the boundary. Since we only have information about the
geometric changes of the boundary, we can also only make assumptions about
geometric changes or displacements in the interior.
As we pointed out in Section~\ref{sec:auxfields}, we borrow the idea from the
Elastic Mesh Update Method (EMUM) to use an elastic material law to determine
the interior displacements.

To keep this description simple, we chose linear elasticity for the smoothing
equation, which was sufficient for all of our test cases. However, more
elaborate equations could be used in its place.
We use the
symbol $\varphi$ for quantities relating to the smoothing equation. The stress
tensor $\boldsymbol{\sigma}^{\varphi}$ that is used in this context is defined
as follows:

\begin{align}
  \boldsymbol{\sigma}^{\varphi}(\vec{r})
  \; &= \; \lambda^{\varphi} (\boldsymbol{\nabla} \cdot \vec{r}) \bm{I}
           + 2 \mu^{\varphi} \, \text{sym} \boldsymbol{\nabla} \vec{r} \, ,
\end{align}

with the virtual Lam\'e parameters $\lambda^{\varphi}$ and $\mu^{\varphi}$.


\subsection{Forward simulation}

We will now formulate the differential equations that are required to run a
forward simulation of thermoelasticity with initial displacements.
We use the solution space $\mathcal{S} = C^2(\Omega)$, with our reference domain $\Omega$.
The spatial dimensionality will be referred to as $d$, such that $\mathcal{S}^d$
is a multi-dimensional solution space.
The momentum conservation is ensured by equation

\begin{alignat}{3}
  \label{eqn:momentum}
  - \nabla \cdot \firstpkref(\vec{r}, \vec{u}, \theta)
  - \bar{\vec{f}} \bar{J}(\vec{u}) \; &= \; \bm{0}
    && \qquad , \, \text{on} \; \Omega \setminus \partial \Omega \, ,
\end{alignat}

with the first Piola-Kirchhoff stress tensor in the reference configuration
$\firstpkref : \mathcal{S}^d \times \mathcal{S}^d \times \mathcal{S} \to
\left[C^1(\Omega)\right]^{d \times d}$ as a function of the
initial displacements $\vec{r} \in \mathcal{S}^d$,
the equilibrium displacements $\vec{u} \in \mathcal{S}^d$, and
the temperature $\theta \in \mathcal{S}$.

$\bar{\vec{f}}$ is a field of forces acting in the equilibrium
configuration, which is defined in the reference configuration.

The Dirichlet and Neumann boundary conditions for this
equation are given by

\begin{alignat}{3}
  \label{eqn:bnddirichlet}
  \vec{u} \; &= \; \vec{u}^D
    && \qquad , \, \text{on} \; \Gamma_D \, , \\
  \label{eqn:bndneumann}
  \firstpkref(\vec{r}, \vec{u}, \theta) \vec{n}
    - \bar{\vec{h}} \bar{J}(\vec{u})
      \left| \bar{\bm{F}}^{-T}(\vec{u}) \vec{n} \right|
      \; &= \; \bm{0}
    && \qquad , \, \text{on} \; \Gamma_N \, ,
\end{alignat}

where $\Gamma_D \cup \Gamma_N = \partial \Omega$, which is the full boundary of
$\Omega$. $\vec{u}^D$ is a prescribed boundary displacement and $\bar{\vec{h}}$
is a force acting on the boundary in the equilibrium configuration.

The smoothing equation for $\vec{r}$ is in this case also a type of momentum
conservation equation, but for the virtual smoothing stress
$\boldsymbol{\sigma}^{\varphi}$:

\begin{alignat}{3}
  \label{eqn:smoothing}
  - \nabla \cdot \boldsymbol{\sigma}^{\varphi}(\vec{r}) \; &= \; \bm{0}
    && \qquad , \, \text{on} \; \Omega \setminus \partial \Omega \, .
\end{alignat}

There are no Neumann boundary conditions for the smoothing equation, as the
boundary shape is completely prescribed. Therefore, with the boundary
displacements $\vec{r}^D$ given in the initial configuration, we have

\begin{alignat}{3}
  \label{eqn:bndsmoothing}
  \vec{r} \; &= \; \vec{r}^D
    && \qquad , \, \text{on} \; \partial \Omega \, .
\end{alignat}

Stationary energy conservation is fulfilled by the heat conduction equation

\begin{alignat}{3}
  \label{eqn:heat}
  - \nabla \cdot \heatfluxref(\vec{r}, \vec{u}, \theta)
  - \bar{g} \bar{J}(\vec{u}) \; &= \; 0
    && \qquad , \, \text{on} \; \Omega \setminus \partial \Omega \, ,
\end{alignat}

with heat flux vector in the reference configuration
$\heatfluxref : \mathcal{S}^d \times \mathcal{S}^d \times \mathcal{S} \to
\left[C^1(\Omega)\right]^d$.
$\bar{g}$ is a heat source field in the
equilibrium configuration.
The Dirichlet and Neumann boundary conditions,

\begin{alignat}{3}
  \label{eqn:bndheatd}
  \theta \; &= \; \theta^D
    && \qquad , \, \text{on} \; \Gamma^{\theta}_D \, , \\
  \label{eqn:bndheatn}
  \heatfluxref(\vec{r}, \vec{u}, \theta) \cdot \vec{n}
    - \bar{k} \bar{J}(\vec{u})
      \left| \bar{\bm{F}}^{-T}(\vec{u}) \vec{n} \right|
      \; &= \; 0
    && \qquad , \, \text{on} \; \Gamma^{\theta}_N \, ,
\end{alignat}

can be used to set either fixed temperature values in $\theta^D$ or a heat flux normal to the
surface in the equilibrium configuration ($\bar{k}$).
Analogous to the boundary conditions of the momentum equation, we have
$\Gamma^{\theta}_D \cup \Gamma^{\theta}_N = \partial \Omega$.

In summary, our system of equations becomes

\begin{alignat}{3}
  \tag{\ref{eqn:momentum}}
  - \nabla \cdot \firstpkref(\vec{r}, \vec{u}, \theta)
  - \bar{\vec{f}} \bar{J}(\vec{u}) \; &= \; \bm{0}
    && \qquad , \, \text{on} \; \Omega \setminus \partial \Omega \, , \\
  \tag{\ref{eqn:smoothing}}
  - \nabla \cdot \boldsymbol{\sigma}^{\varphi}(\vec{r}) \; &= \; \bm{0}
    && \qquad , \, \text{on} \; \Omega \setminus \partial \Omega \, , \\
  \tag{\ref{eqn:heat}}
  - \nabla \cdot \heatfluxref(\vec{r}, \vec{u}, \theta)
  - \bar{g} \bar{J}(\vec{u}) \; &= \; 0
    && \qquad , \, \text{on} \; \Omega \setminus \partial \Omega \, , \\
  \tag{\ref{eqn:bnddirichlet}}
  \vec{u} \; &= \; \vec{u}^D
    && \qquad , \, \text{on} \; \Gamma_D \, , \\
  \tag{\ref{eqn:bndneumann}}
  \firstpkref(\vec{r}, \vec{u}, \theta) \vec{n}
    - \bar{\vec{h}} \bar{J}(\vec{u})
      \left| \bar{\bm{F}}^{-T}(\vec{u}) \vec{n} \right|
      \; &= \; \bm{0}
    && \qquad , \, \text{on} \; \Gamma_N \, , \\
  \tag{\ref{eqn:bndsmoothing}}
  \vec{r} \; &= \; \vec{r}^D
    && \qquad , \, \text{on} \; \partial \Omega \, , \\
  \tag{\ref{eqn:bndheatd}}
  \theta \; &= \; \theta^D
    && \qquad , \, \text{on} \; \Gamma^{\theta}_D \, , \\
  \tag{\ref{eqn:bndheatn}}
  \heatfluxref(\vec{r}, \vec{u}, \theta) \cdot \vec{n}
    - \bar{k} \bar{J}(\vec{u})
      \left| \bar{\bm{F}}^{-T}(\vec{u}) \vec{n} \right|
      \; &= \; 0
    && \qquad , \, \text{on} \; \Gamma^{\theta}_N \, .
\end{alignat}

\subsection{Inverse simulation}
\label{sec:invsim}

In order to run an inverse simulation instead of the forward simulation, very
few things need to be changed. Notably, the differential equations for the
conservation laws, i.e., Equations~(\ref{eqn:momentum}), (\ref{eqn:smoothing}),
and (\ref{eqn:heat}),
stay in place.
The boundary conditions for the heat equation,
Equations~(\ref{eqn:bndheatd})~and~(\ref{eqn:bndheatn}), also remain unchanged.
In order to switch from prescribing the boundary in the initial configuration to
prescribing it in the equilibrium configuration,
Equations~(\ref{eqn:bnddirichlet}), (\ref{eqn:bndneumann}),
and (\ref{eqn:bndsmoothing}) are replaced by the new boundary conditions

\begin{alignat}{3}
  \label{eqn:bnddirichletinv}
  \vec{u} \; &= \; \vec{u}^D
    && \qquad , \, \text{on} \; \partial \Omega \, , \\
  \label{eqn:bndneumanninv}
  \firstpkref(\vec{r}, \vec{u}, \theta) \vec{n}
    - \bar{\vec{h}} \bar{J}(\vec{u})
      \left| \bar{\bm{F}}^{-T}(\vec{u}) \vec{n} \right|
      \; &= \; \bm{0}
    && \qquad , \, \text{on} \; \partial \Omega \, .
\end{alignat}

This specifies both Dirichlet and Neumann boundary conditions for the
elasticity equations on the whole boundary. This over-specification
makes it possible to under-specify boundary conditions for the smoothing
equation.
Notably, both $\vec{u}^D$ and $\bar{\vec{h}}$ now have to be available on the
full boundary $\partial \Omega$. If the shape $\vec{u}^D$ is known for an unconstrained state,
$\bar{\vec{h}}$ is zero on the full boundary and
Equation~(\ref{eqn:bndneumanninv}) simplifies to
$\firstpkref(\vec{r}, \vec{u}, \theta) \vec{n} = 0$.

The full system of equations for the inverse formulation reads

\begin{alignat}{3}
  \tag{\ref{eqn:momentum}}
  - \nabla \cdot \firstpkref(\vec{r}, \vec{u}, \theta)
  - \bar{\vec{f}} \bar{J}(\vec{u}) \; &= \; \bm{0}
    && \qquad , \, \text{on} \; \Omega \setminus \partial \Omega \, , \\
  \tag{\ref{eqn:smoothing}}
  - \nabla \cdot \boldsymbol{\sigma}^{\varphi}(\vec{r}) \; &= \; \bm{0}
    && \qquad , \, \text{on} \; \Omega \setminus \partial \Omega \, , \\
  \tag{\ref{eqn:heat}}
  - \nabla \cdot \heatfluxref(\vec{r}, \vec{u}, \theta)
  - \bar{g} \bar{J}(\vec{u}) \; &= \; 0
    && \qquad , \, \text{on} \; \Omega \setminus \partial \Omega \, , \\
  \tag{\ref{eqn:bnddirichletinv}}
  \vec{u} \; &= \; \vec{u}^D
    && \qquad , \, \text{on} \; \partial \Omega \, , \\
  \tag{\ref{eqn:bndneumanninv}}
  \firstpkref(\vec{r}, \vec{u}, \theta) \vec{n}
    - \bar{\vec{h}} \bar{J}(\vec{u})
      \left| \bar{\bm{F}}^{-T}(\vec{u}) \vec{n} \right|
      \; &= \; \bm{0}
    && \qquad , \, \text{on} \; \partial \Omega \, , \\
  \tag{\ref{eqn:bndheatd}}
  \theta \; &= \; \theta^D
    && \qquad , \, \text{on} \; \Gamma^{\theta}_D \, , \\
  \tag{\ref{eqn:bndheatn}}
  \heatfluxref(\vec{r}, \vec{u}, \theta) \cdot \vec{n}
    - \bar{k} \bar{J}(\vec{u})
      \left| \bar{\bm{F}}^{-T}(\vec{u}) \vec{n} \right|
      \; &= \; 0
    && \qquad , \, \text{on} \; \Gamma^{\theta}_N \, .
\end{alignat}

\subsection{Discretized Form of Equations}

We solve the equations using the finite element method. As a preparation, we need
to find suitable interpolation spaces and formulate the discretized weak form of
the equations. We will sketch the interpolation spaces for both the standard
finite element method with linear basis functions as well as isogeometric
analysis with B-spline or NURBS basis functions. Other spaces can also be used but have not
been tested.

We will first construct the interpolation space for the finite element method
with simplex P1 (linear) basis functions. We use a reference element defined by
a compact set $\hat{\Omega} \subset \mathbb{R}^d$.
For the individual elements $e \in E$, where $E$ is the discrete set of all element
indices, we define a linear projector into physical space

\begin{align}
  &\mathcal{P}_e : \hat{\Omega} \to \Omega_e \, , \\
  &\mathcal{P}_e \in \left[ \mathbb{P}_1(\hat{\Omega}) \right]^d \, ,
\end{align}

with the space of linear polynomials $\mathbb{P}_1$,
where $\Omega_e \subseteq \Omega$, and especially

\begin{align}
   \lambda\left(\Omega_e \cap \Omega_f\right) \; &= \; 0 \, , \quad
     \forall e, f \in E : e \ne f \, , \\
  \bigcup_{e \in E} \bar{\Omega}_e \; &= \; \Omega \, ,
\end{align}

where $\lambda(\square)$ is the Lebesgue measure.
This means that the physical elements $\Omega_e$ do not overlap except on their
boundaries,
and when they are combined, they make up the full domain $\Omega$.
The finite element interpolation space can be constructed using the projectors
as

\begin{align}
  \mathcal{I}_{\text{lin}} \; = \; \left\{ f \in C(\Omega) \,
  \middle| \, f|_{\Omega_e} \circ \mathcal{P}_e \in \mathbb{P}_1(\hat{\Omega}),
              \forall e \in E \right\} \, ,
\end{align}

which makes it a space of piecewise linear continuous functions on $\Omega$.

As an alternative, we define the spline interpolation space for isogeometric
analysis as

\begin{align}
  \mathcal{I}_{\text{spline}} \; &= \; \text{span} \, B_i \, ,
\end{align}

with a compact reference space $\tilde{\Omega} \subset \mathbb{R}^d$ and the
spline basis functions $B_i : \tilde{\Omega} \to \mathbb{R}$, which could be,
e.g., B-Spline or NURBS basis functions \cite{Hughes2005}.

In the following, we will refer to the interpolation space simply as
$\mathcal{I}$, assuming that one of the presented options is used.
Based on this we define the solution spaces

\begin{align}
  \mathcal{S}_{\vec{u}} \; &= \; \left\{ \vec{v} \in \mathcal{I}^d \; \middle| \;
                        \vec{v}|_{\Gamma_D} = \vec{u}^D \right\} \, , \\
  \mathcal{S}_{\vec{r}} \; &= \; \left\{ \vec{v} \in \mathcal{I}^d \; \middle| \;
                        \vec{v}|_{\partial \Omega} = \vec{r}^D \right\} \, , \\
  \mathcal{S}_{\theta} \; &= \; \left\{ \vec{s} \in \mathcal{I} \; \middle| \;
                        s|_{\Gamma^{\theta}_D} = \theta^D \right\} \, ,
\end{align}

and the test spaces

\begin{align}
  \mathcal{T}_{\vec{U}} \; &= \; \left\{ \vec{w} \in \mathcal{I}^d \; \middle| \;
                         \vec{w}|_{\Gamma_D} = \vec{0} \right\} \, , \\
  \mathcal{T}_{\vec{R}} \; &= \; \left\{ \vec{w} \in \mathcal{I}^d \; \middle| \;
                         \vec{w}|_{\partial \Omega} = \vec{0} \right\} \, , \\
  \mathcal{T}_{\Theta} \; &= \; \left\{ \vec{s} \in \mathcal{I} \; \middle| \;
                        s|_{\Gamma^{\theta}_D} = 0 \right\} \, .
\end{align}

For simulations in forward mode, we obtain the following discretized weak
formulation:
\textit{Find $\vec{u}^h \in \mathcal{S}_{\vec{u}}$,
$\vec{r}^h \in \mathcal{S}_{\vec{r}}$, and $\theta^h \in \mathcal{S}_{\theta}$,
such that}

\begin{equation}
  \begin{aligned}
  \label{eqn:discweakforward}
  0 \; =& \; \int_{\Omega} \left[ \boldsymbol{\nabla} \vec{U}^h
               : \firstpkref(\vec{r}^h, \vec{u}^h, \theta^h)
               - \left(\vec{U}^h \cdot \bar{\vec{f}}\right) \bar{J}(\vec{u}^h)
               \right] \, \text{d}\Omega
             - \int_{\Gamma_N} \left(\vec{U}^h \cdot \bar{\vec{h}}\right) \bar{J}(\vec{u}^h)
               \left| \bar{\bm{F}}^{-T}(\vec{u}^h) \vec{n} \right|
               \, \text{d}\Gamma \\
        & \; + \int_{\Omega} \boldsymbol{\nabla} \vec{R}^h
               : \boldsymbol{\sigma}^{\varphi}(\vec{r}^h)
               \, \text{d}\Omega \\
        & \; + \int_{\Omega} \left[ \boldsymbol{\nabla} \Theta^h
               \cdot \heatfluxref(\vec{r}^h, \vec{u}^h, \theta^h)
               - \Theta^h \bar{g} \bar{J}(\vec{u}^h)
               \right] \, \text{d}\Omega
             - \int_{\Gamma^{\theta}_N} \Theta^h \bar{k} \bar{J}(\vec{u}^h)
               \left| \bar{\bm{F}}^{-T}(\vec{u}^h) \vec{n} \right|
               \, \text{d}\Gamma \, ,
  \end{aligned}
\end{equation}

\textit{for all $\vec{U}^h \in \mathcal{T}_{\vec{U}}$,
$\vec{R}^h \in \mathcal{T}_{\vec{R}}$, and
$\Theta^h \in \mathcal{T}_{\Theta}$.}

For the inverse simulation, we have to make changes to the solution spaces for the
displacement fields:

\begin{align}
  \bar{\mathcal{S}}_{\vec{u}} \; &= \; \left\{ \vec{v} \in \mathcal{I}^d \; \middle| \;
                        \vec{v}|_{\partial \Omega} = \vec{u}^D \right\} \, , \\
  \bar{\mathcal{S}}_{\vec{r}} \; &= \; \mathcal{I}^d \, .
\end{align}

Additionally, the test space for the momentum conservation equation has to be
adjusted to

\begin{align}
  \bar{\mathcal{T}}_{\vec{U}} \; &= \; \mathcal{I}^d \, ,
\end{align}

to incorporate the full domain including its boundary.

The discretized weak formulation of the inverse problem is defined as:
\textit{Find $\vec{u}^h \in \bar{\mathcal{S}}_{\vec{u}}$,
$\vec{r}^h \in \bar{\mathcal{S}}_{\vec{r}}$, and $\theta^h \in \mathcal{S}_{\theta}$,
such that}

\begin{equation}
  \begin{aligned}
  \label{eqn:discweakinverse}
  0 \; =& \; \int_{\Omega} \left[ \boldsymbol{\nabla} \vec{U}^h
               : \firstpkref(\vec{r}^h, \vec{u}^h, \theta^h)
               - \left(\vec{U}^h \cdot \bar{\vec{f}}\right) \bar{J}(\vec{u}^h)
               \right] \, \text{d}\Omega
             - \int_{\partial \Omega} \left(\vec{U}^h \cdot \bar{\vec{h}}\right) \bar{J}(\vec{u}^h)
               \left| \bar{\bm{F}}^{-T}(\vec{u}^h) \vec{n} \right|
               \, \text{d}\Gamma \\
        & \; + \int_{\Omega} \boldsymbol{\nabla} \vec{R}^h
               : \boldsymbol{\sigma}^{\varphi}(\vec{r}^h)
               \, \text{d}\Omega \\
        & \; + \int_{\Omega} \left[ \boldsymbol{\nabla} \Theta^h
               \cdot \heatfluxref(\vec{r}^h, \vec{u}^h, \theta^h)
               - \Theta^h \bar{g} \bar{J}(\vec{u}^h)
               \right] \, \text{d}\Omega
             - \int_{\Gamma^{\theta}_N} \Theta^h \bar{k} \bar{J}(\vec{u}^h)
               \left| \bar{\bm{F}}^{-T}(\vec{u}^h) \vec{n} \right|
               \, \text{d}\Gamma
  \end{aligned}
\end{equation}

\textit{for all $\vec{U}^h \in \bar{\mathcal{T}}_{\vec{U}}$,
$\vec{R}^h \in \mathcal{T}_{\vec{R}}$, and
$\Theta^h \in \mathcal{T}_{\Theta}$.}

When compared to the equation for the forward simulation in
(\ref{eqn:discweakforward}), Equation~(\ref{eqn:discweakinverse}) seems almost
identical. The only change is in the integral limits of the boundary term for
the momentum conservation equation, which now contain the full boundary
$\partial \Omega$. In
cases where $\bar{h} \equiv 0$, this term is removed and the equations become
identical.

\section{Derivation of Constitutive Equations}

So far, we have left the first Piola-Kirchhoff stress tensor $\firstpkref$ and the
heat flux $\heatfluxref$ in the reference configuration undefined. We will provide definitions for both of these
quantities in the following sections.

We have already provided some definitions of deformation measures in
Figure~\ref{fig:configurations}. In addition, we will now introduce the
deformation gradient tensor $\deftentot$, which describes the full deformation in
the equilibrium configuration with respect to the stress-free configuration:

\begin{align}
  \deftentot \; &:= \; \D{\bar{\vec{x}}}{\hat{\vec{x}}} \, . \\
\end{align}

Analogously, we define $J := \text{det} \, \deftentot$.

Due to the multiple coordinate systems that we use, we also have to define
multiple stress measures. Table~\ref{tab:stressmeasures} gives an overview of
the stress tensors that we use. The tensors $\firstpktot$ and $\secondpktot$
correspond to the first and second Piola-Kirchhoff stress tensors, respectively,
as they are commonly used.

\begin{table}[htbp]
  \centering
  \caption{Overview of stress measures}
  \begin{tabular}{| c || c | c |}
    \hline
    Symbol & relates forces in & to surface elements in \\
    \hline
    \hline
    $\boldsymbol{\sigma}$ & equilibrium configuration $\bar{\vec{x}}$ & equilibrium configuration $\bar{\vec{x}}$ \\
    \hline
    $\firstpktot$ & equilibrium configuration $\bar{\vec{x}}$ & stress-free configuration $\hat{\vec{x}}$ \\
    \hline
    $\secondpktot$ & stress-free configuration $\hat{\vec{x}}$ & stress-free configuration $\hat{\vec{x}}$ \\
    \hline
    $\firstpkref$ & equilibrium configuration $\bar{\vec{x}}$ & reference configuration $\vec{x}$ \\
    \hline
    $\secondpkref$ & reference configuration $\vec{x}$ & reference configuration $\vec{x}$ \\
    \hline
  \end{tabular}
  \label{tab:stressmeasures}
\end{table}

\subsection{Constitutive Laws for the First Piola-Kirchhoff Stress}

In many cases where nonlinear elasticity laws are used, the reference
configuration coincides with the stress-free configuration, i.e., $\firstpktot =
\firstpkref$, as well as, $\secondpktot = \secondpkref$. For us, this is not
the case. Therefore, when we use constitutive laws that are defined with a
stress-free configuration used as reference, we need to carry out a coordinate
transformation to the reference configuration.

\subsubsection{Based on St.\ Venant-Kirchhoff}

The following is a simple constitutive law for thermoelasticity that is based on the
St.\ Venant-Kirchhoff material:

\begin{align}
  \secondpktot \; &= \; \lambda (\mathrm{tr} \, \bm{E}) \bm{I} + 2 \mu \bm{E}
                   - \alpha (\theta - \theta_0) \bm{I} \, ,
\end{align}

where the Lam\'e parameters $\lambda$ and $\mu$ control the elastic behavior and
$\alpha$ controls the thermal expansion (cf. \cite{Surana2015} for a derivation
of this law).
It uses a simple linear thermal expansion law. Such a law
should only be used with small temperature variations or for demonstration
purposes.
For this constitutive law, we can derive the following expression in
$\firstpkref$, which can be inserted into the presented differential
equations\footnote{For a full derivation of this equation, see~\ref{sec:reformstvk}.}:

\begin{align}
  \firstpkref \; &= \;
      \tilde{J} \mathring{J}^{-1}
      \left( \frac{\lambda}2 \mathrm{tr} (\bm{B} - \bm{I}) \bm{I} + \mu (\bm{B} - \bm{I})
    - \alpha (\theta - \theta_0) \bm{I} \right)
      \bar{\bm{F}} \tilde{\bm{F}}^{-1} \mathring{\bm{B}}
             \tilde{\bm{F}}^{-T} \, ,
\end{align}

with $\bm{B} = \bar{\bm{F}} \tilde{\bm{F}}^{-1} \mathring{\bm{B}} \tilde{\bm{F}}^{-T} \bar{\bm{F}}^T$.

\subsubsection{Based on Neo-Hooke}

We also use a constitutive law that is based on a Neo-Hooke material (cf., e.g.,
\cite{Haupt2002}), with
additional thermal expansion terms. This is given as an expression that relates
the Cauchy stress and the left Cauchy-Green deformation tensor $\bm{B}$:

\begin{align}
  J \sigma \; &= \; 2 D_1 J (J-1) \bm{I} + 2 C_1 J^{-\frac23} (\text{dev} \, \bm{B})
                   - \alpha (\theta - \theta_0) \bm{B} \, .
\end{align}

The elastic behavior is controlled through the material parameters $D_1$ and
$C_1$, and the thermal expansion behavior through the parameter $\alpha$.
The thermal expansion term has been chosen this way for consistency with the law
based on St. Venant-Kirchhoff that was presented in the previous section.
For the stress tensor $\firstpkref$ that we require for our differential
equations, we obtain the following
expression\footnote{For a full derivation of this equation, see~\ref{sec:reformneohooke}.}:

\begin{align}
  \firstpkref
  \; &= \; \tilde{J} \mathring{J}^{-1} \left(2 D_1 J (J-1) \bm{I} + 2 C_1 J^{-\frac23} (\text{dev} \, \bm{B})
                 - \alpha (\theta - \theta_0) \bm{B}\right) \bar{\bm{F}}^{-T} \, ,
\end{align}

with $J = \bar{J} \tilde{J}^{-1} \mathring{J}$ and
$\bm{B} = \bar{\bm{F}} \tilde{\bm{F}}^{-1} \mathring{\bm{B}} \tilde{\bm{F}}^{-T} \bar{\bm{F}}^T$.

\subsection{Constitutive Law for the Heat Flux}

We have formulated the differential equations for an arbitrary heat flux vector
$\heatfluxref(\vec{r}, \vec{u}, \theta)$. We have done so to point out that multiple
options of defining the heat flux are available.
For instance, the heat flux could be chosen to be consistent with regular
formulations where the heat flux is a linear function of the temperature
gradient with respect to the stress-free or equilibrium coordinate system.
However, the only definition that we have used is one where the heat flux is
linear with respect to the temperature gradient in the reference configuration:

\begin{align}
  \heatfluxref(\vec{r}, \vec{u}, \theta) \; &= \; - \kappa \nabla \theta \, ,
\end{align}

with scalar thermal conductivity $\kappa$.

\section{Numerical Examples}

We describe two different test cases in this section. The first test case,
described in Section~\ref{sec:elasticbeam}, is meant to be a simple
demonstration of the capabilities of the inverse formulation when it is applied
to isothermal elasticity. In the second test
case, Section~\ref{sec:thermobody}, we employ all the concepts that were
introduced in the previous chapters to solve a thermoelastic problem.

\subsection{Elastic beam under gravitational forces}
\label{sec:elasticbeam}

In this numerical example, we answer the question of how a clamped elastic beam
would have to be shaped, such that it takes on a completely straight shape
when subjected to gravity.
The test case demonstrates the simplest
manner how the inverse formulation could be used, since temperature is ignored
and initial stresses are not used.

\begin{figure}[htbp]
  \centering
  \tiny
  \begin{tabular}{c c c}
  & $\vec{r} = 0, \; \boldsymbol{\sigma} \cdot \vec{n} = 0$ \vspace{.3cm} & \\
  $\begin{aligned} \vec{r} &= 0 \\ \vec{u} &= 0 \\ \\ \\ \end{aligned}$ &
    \includegraphics[width=12cm]{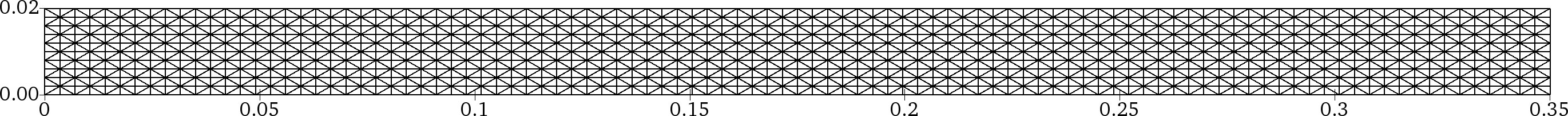} &
    $\begin{aligned} \vec{r} &= 0 \\ \boldsymbol{\sigma} \cdot \vec{n} &= 0 \\ \\ \\ \end{aligned}$ \\
  & $\vec{r} = 0, \; \boldsymbol{\sigma} \cdot \vec{n} = 0$ & \\
  \end{tabular}
  \caption{Simulation mesh for an elastic beam of size $0.02 \times 0.35$ in the reference
  configuration. The boundary conditions are shown for the forward simulation,
  where the beam is attached to a wall on the left-hand side, and free to move
  everywhere else.}
  \label{fig:gravity.setup}
\end{figure}

We consider a simplified two-dimensional beam. The simulation mesh in the reference
configuration, along with boundary conditions for the forward simulation,
is shown in Figure~\ref{fig:gravity.setup}.

We use the Dirichlet boundary condition $\vec{u} =
0$ on the left-hand side boundary for the equilibrium displacement field. This
makes sure that the beam stays attached to the wall.
Additionally, we use the Dirichlet boundary condition $\vec{r} = 0$ on all
boundaries for the initial displacement field.
This implies that the smoothing equation will ensure that the initial and
reference configurations are identical\footnote{The
initial displacement field as well as the smoothing equation are not necessary
in the forward simulation, but we use them anyway to highlight the differences
between forward and inverse setups.}.
The values for the material properties and the gravitational force are shown in
Table~\ref{tab:gravity.simparams}.

\begin{table}[htbp]
  \centering
  \caption{Simulation parameters}
  \begin{tabular}{| l | r |}
    \hline
    Elasticity model & St. Venant-Kirchhoff \\
    \hline
    $\lambda$ & $2 \times 10^6$ \\
    \hline
    $\mu$ & $0.5 \times 10^6$ \\
    \hline
    $\bar{\vec{f}} = \rho \vec{g}$ & $(0, -2000)^T$ \\
    \hline
  \end{tabular}
  \label{tab:gravity.simparams}
\end{table}

\begin{figure}[htbp]
  \centering
  \includegraphics[width=12cm]{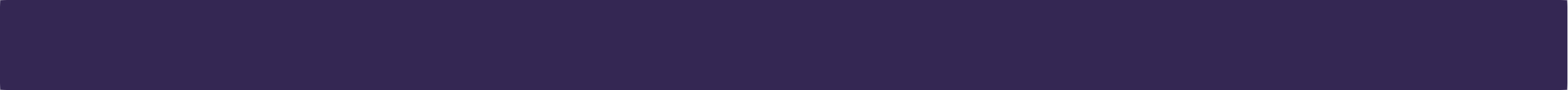}

  \vspace{1cm}

  \includegraphics[width=12cm]{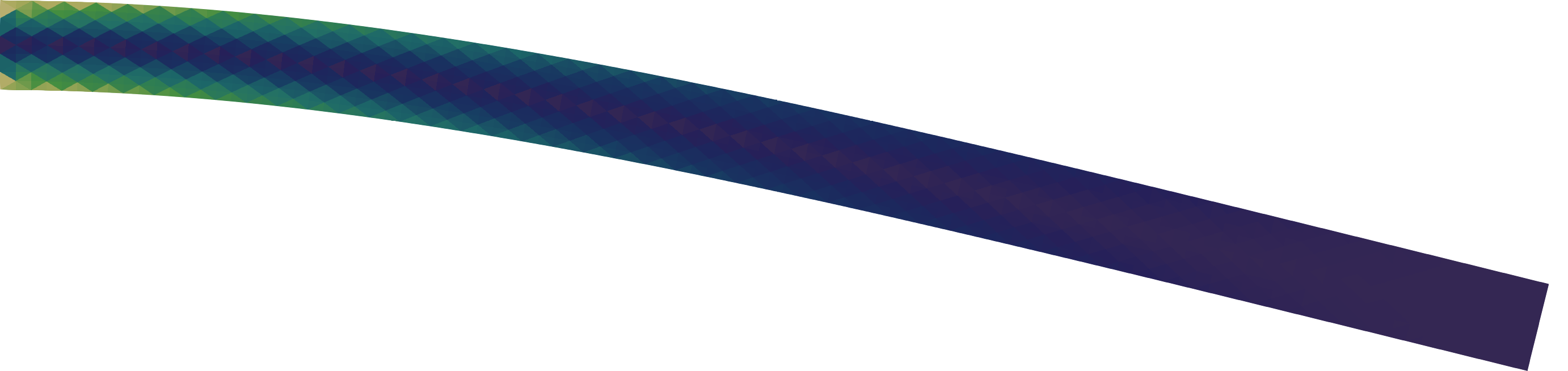}

  \vspace{1cm}

  \includegraphics[width=7cm]{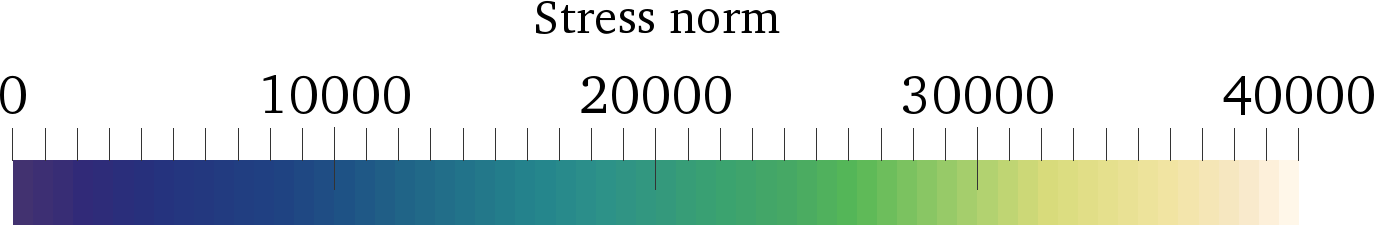}
  \caption{Spectral norm of initial stress tensor in initial configuration (top)
  and residual (equilibrium) stress tensor in equilibrium configuration (bottom)
  for elastic beam in forward simulation. The lower
  right-hand corner of the beam is displaced by $(-0.00895, -0.06279)^T$.}
  \label{fig:gravity.forward.snorm}
\end{figure}

In the equilibrium configuration, the lower right-hand corner of the beam is
displaced by $-8.95 \times 10^{-3}$ in horizontal direction and
$-62.79 \times 10^{-3}$ in vertical
direction. Figure~\ref{fig:gravity.forward.snorm} shows the stresses in the
beam in the initial and the equilibrium configurations.

\begin{figure}[htbp]
  \centering
  \tiny
  \begin{tabular}{c c c}
  & $\boldsymbol{\sigma} \cdot \vec{n} = 0, \; \vec{u} = 0$ \vspace{.3cm} & \\
  $\begin{aligned} \vec{r} &= 0 \\ \vec{u} &= 0 \\ \\ \\ \end{aligned}$ &
    \includegraphics[width=12cm]{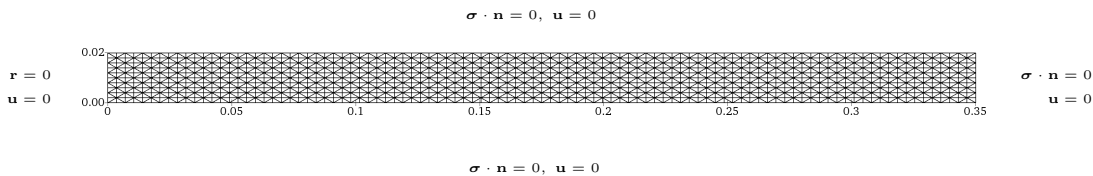} &
    $\begin{aligned} \boldsymbol{\sigma} \cdot \vec{n} &= 0 \\ \vec{u} &= 0 \\ \\ \\ \end{aligned}$ \\
  & $\boldsymbol{\sigma} \cdot \vec{n} = 0, \; \vec{u} = 0$ & \\
  \end{tabular}
  \caption{For the inverse simulation, the boundary conditions are changed, such
  that we can solve for the equilibrium displacement field $\vec{u}$.}
  \label{fig:gravity.setup.inverse}
\end{figure}

For the inverse problem, we remove most of the Dirichlet boundary conditions on
the initial displacement field $\vec{r}$. We only keep the boundary condition of
$\vec{r} = 0$ on the left-hand side boundary. This makes sure that the shape can
be attached to a wall. We use Dirichlet boundary conditions for the equilibrium
displacement field on all boundaries to completely prescribe the equilibrium
shape. The mesh and boundary conditions are shown in
Figure~\ref{fig:gravity.setup.inverse}.

\begin{figure}[htbp]
  \centering
  \includegraphics[width=12cm]{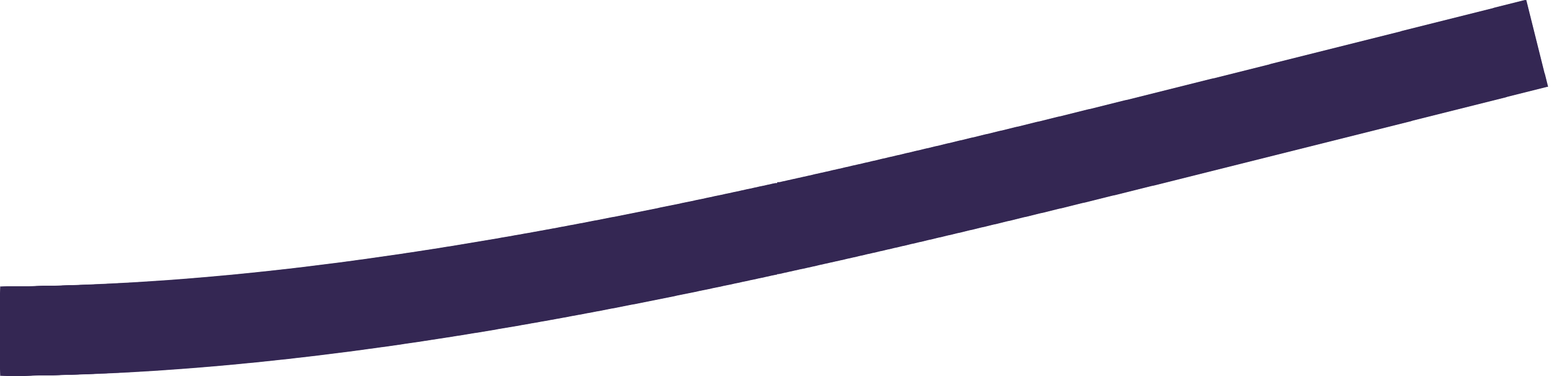}

  \vspace{1cm}

  \includegraphics[width=12cm]{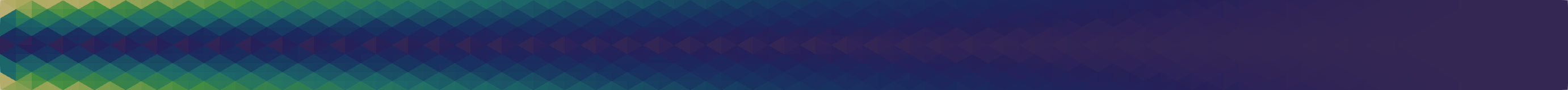}

  \vspace{1cm}

  \includegraphics[width=7cm]{fig-gravity-legend-snorm.png}
  \caption{Spectral norm of initial Cauchy stress tensor in initial configuration (top)
  and residual (equilibrium) Cauchy stress tensor in equilibrium configuration (bottom)
  for elastic beam in inverse simulation. The upper
  right-hand corner of the beam is displaced by $(-0.00923, 0.06397)^T$.}
  \label{fig:gravity.inverse.snorm}
\end{figure}

The solution to the inverse problem is shown in
Figure~\ref{fig:gravity.inverse.snorm}. The
upper right-hand corner of the beam in the initial configuration is displaced by
$-9.23 \times 10^{-3}$ in horizontal direction and
$63.97 \times 10^{-3}$ in vertical direction. While the shape may appear to be
a mirror image of the solution to the forward problem
(Figure~\ref{fig:gravity.forward.snorm}, bottom),
the numbers prove that this is not the case.

\subsection{Thermoelastic body under thermal stresses}
\label{sec:thermobody}

\begin{figure}[htbp]
  \centering
  \includegraphics[width=6cm]{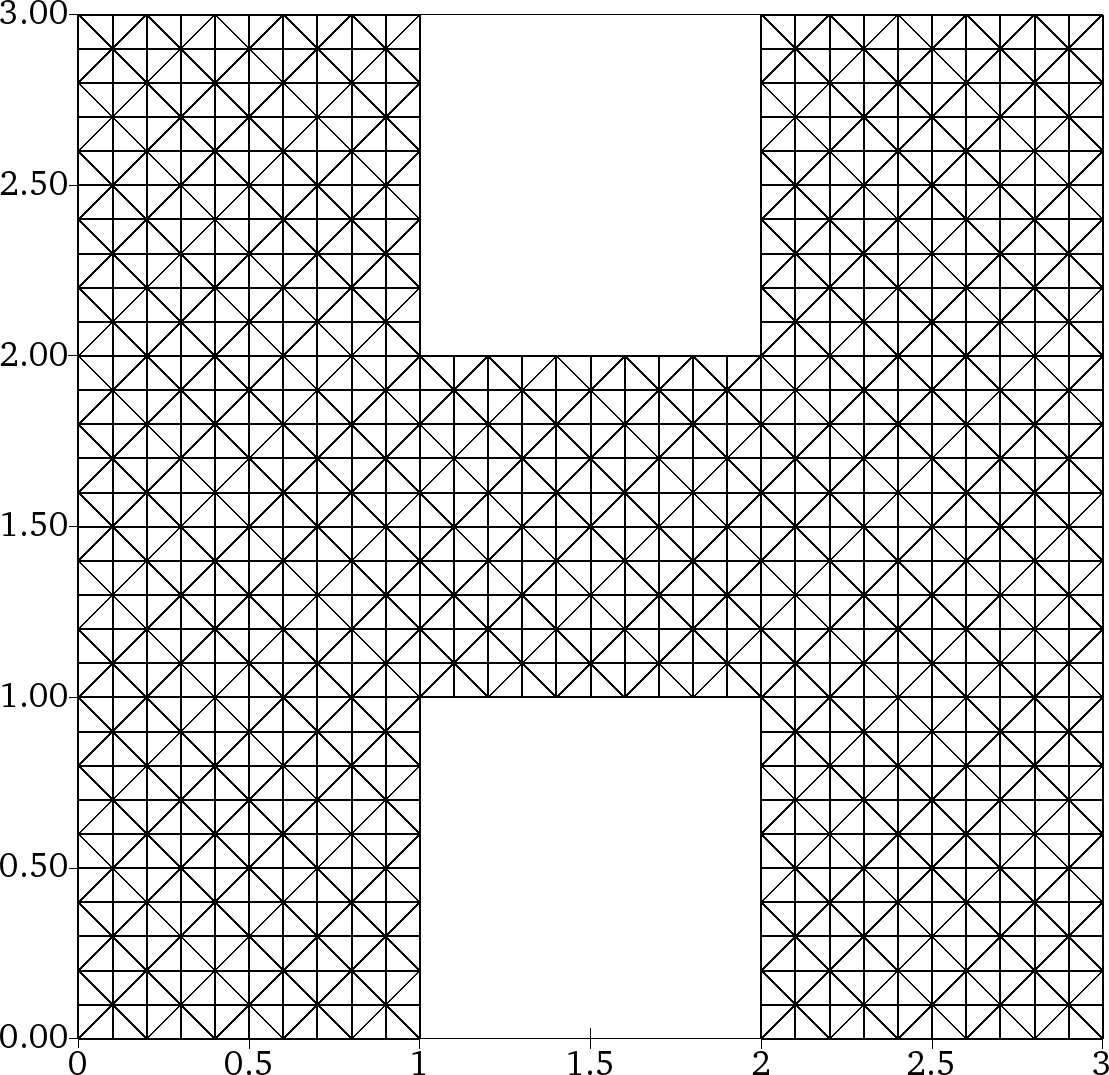}
  \caption{Simulation mesh for thermoelastic body in reference configuration.}
  \label{fig:thermal.setup}
\end{figure}

Our principal test case involves stresses caused by thermal expansion or
contraction. The idea is that we start with an inhomogeneously heated body (such
as an injection molding part that has been partially cooled from the outside) that is
cooled down to an equilibrium temperature. Our example geometry is shown in
Figure~\ref{fig:thermal.setup}.

\begin{figure}[htbp]
  \centering
  \includegraphics[height=6cm]{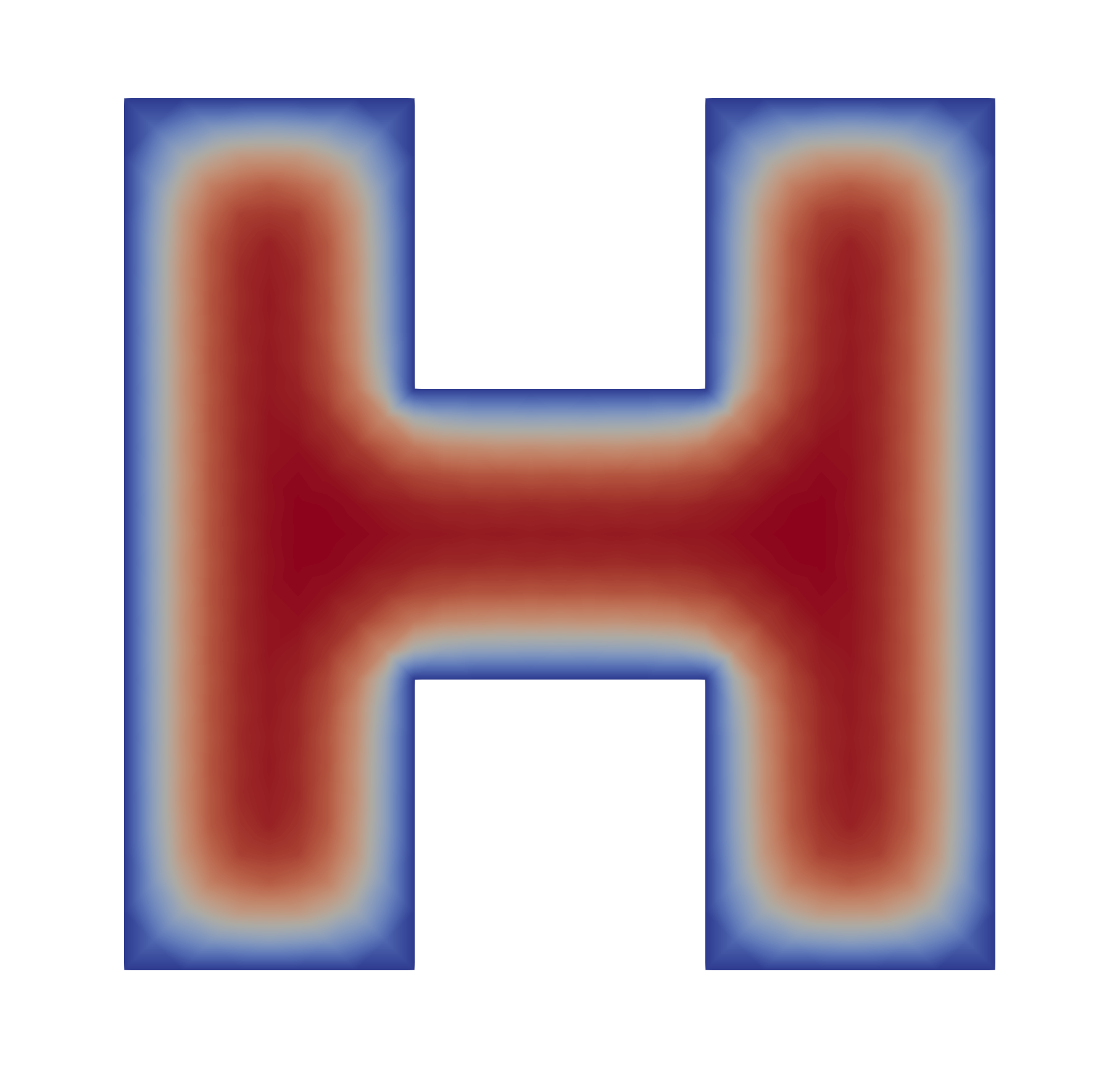}
  \hspace{.4cm}
  \includegraphics[height=6cm]{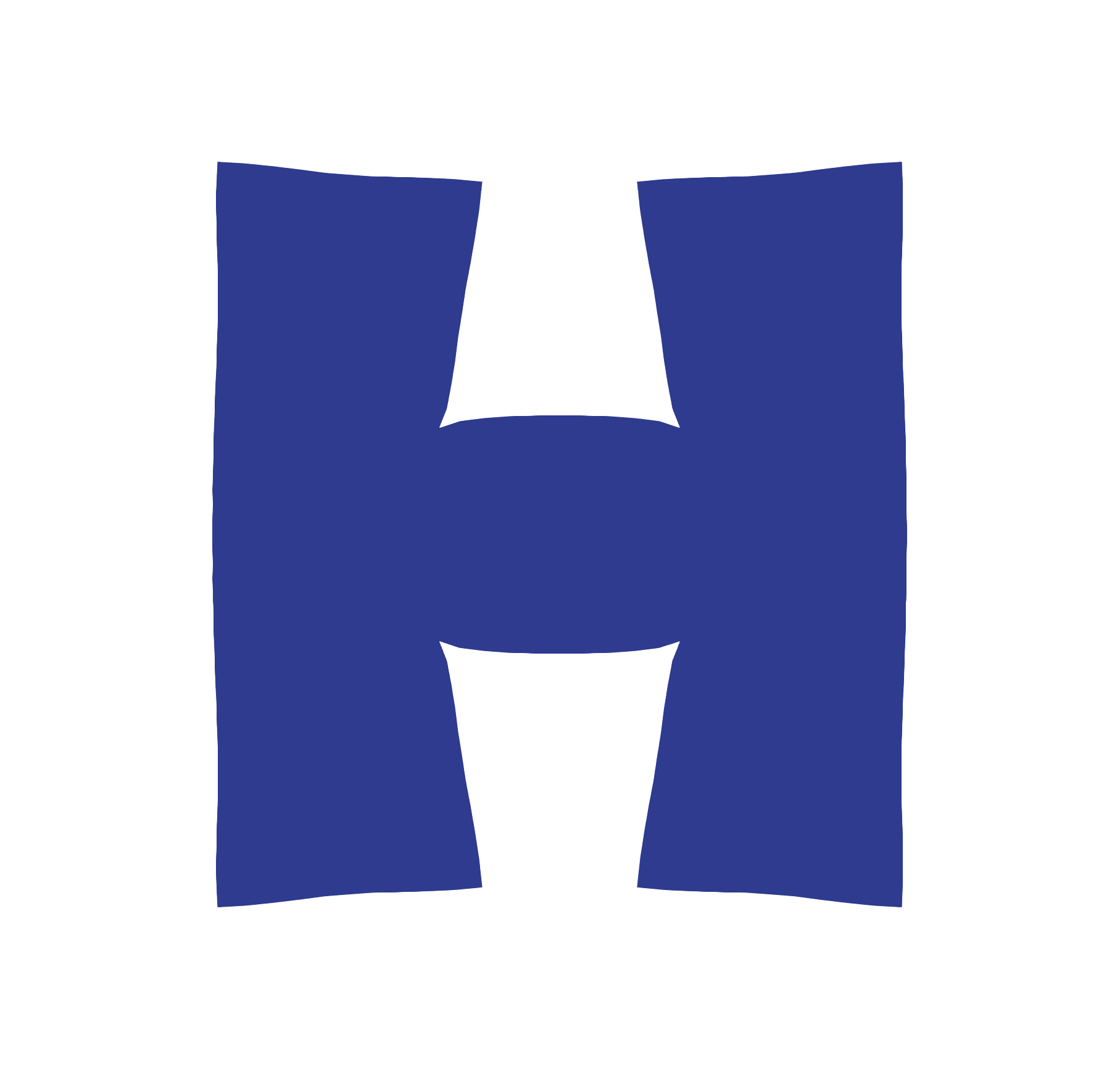}
  \hspace{.4cm}
  \includegraphics[height=6cm]{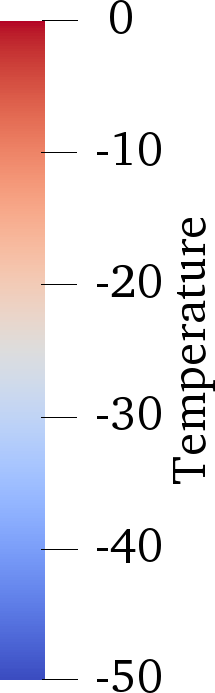}
  \caption{Initial temperature in initial configuration (left)
  and equilibrium temperature in equilibrium configuration (right)
  for thermoelastic body in forward simulation.}
  \label{fig:forward.theta}
\end{figure}

We created an artificial temperature field as the initial condition of our
simulation.
To do so we solved the standard
heat conduction equation, configured
according to Table~\ref{tab:heatsimparams}. Using Dirichlet boundary conditions,
a hot body is slowly cooled down.
Since we required an inhomogeneous temperature field, this simulation was
stopped after five time steps. The result of this simulation is shown on the
left-hand side in Figure~\ref{fig:forward.theta}.

\begin{table}[htbp]
  \centering
  \caption{Heat equation parameters}
  \begin{tabular}{| l | r |}
    \hline
    $\kappa$ & $0.41$ \\
    \hline
    Initial temperature & $0$ \\
    \hline
    Boundary temperature & $-50$ \\
    \hline
    $\Delta t$ & $0.01$ \\
    \hline
    Number of time steps & $5$ \\
    \hline
  \end{tabular}
  \label{tab:heatsimparams}
\end{table}

Table~\ref{tab:simparams} shows the material parameters used for the
thermoelastic forward simulation. One should note that these values do not
correspond to any existing material but are chosen for qualitative testing
purposes.

\begin{table}[htbp]
  \centering
  \caption{Elasticity parameters}
  \begin{tabular}{| l | r |}
    \hline
    Elasticity model & St. Venant-Kirchhoff \\
    \hline
    $\lambda$ & $0.01$ \\
    \hline
    $\mu$ & $100.0$ \\
    \hline
    $\alpha$ & $1.0$ \\
    \hline
    $\kappa$ & $0.2$ \\
    \hline
  \end{tabular}
  \label{tab:simparams}
\end{table}

Starting from the result of the heat conduction simulation, the actual
thermoelastic forward simulation can be run. Figure~\ref{fig:forward.theta} shows
the change that the shape undergoes when the body cools down to a homogeneous
temperature.

\begin{figure}[htbp]
  \centering
  \includegraphics[height=6cm]{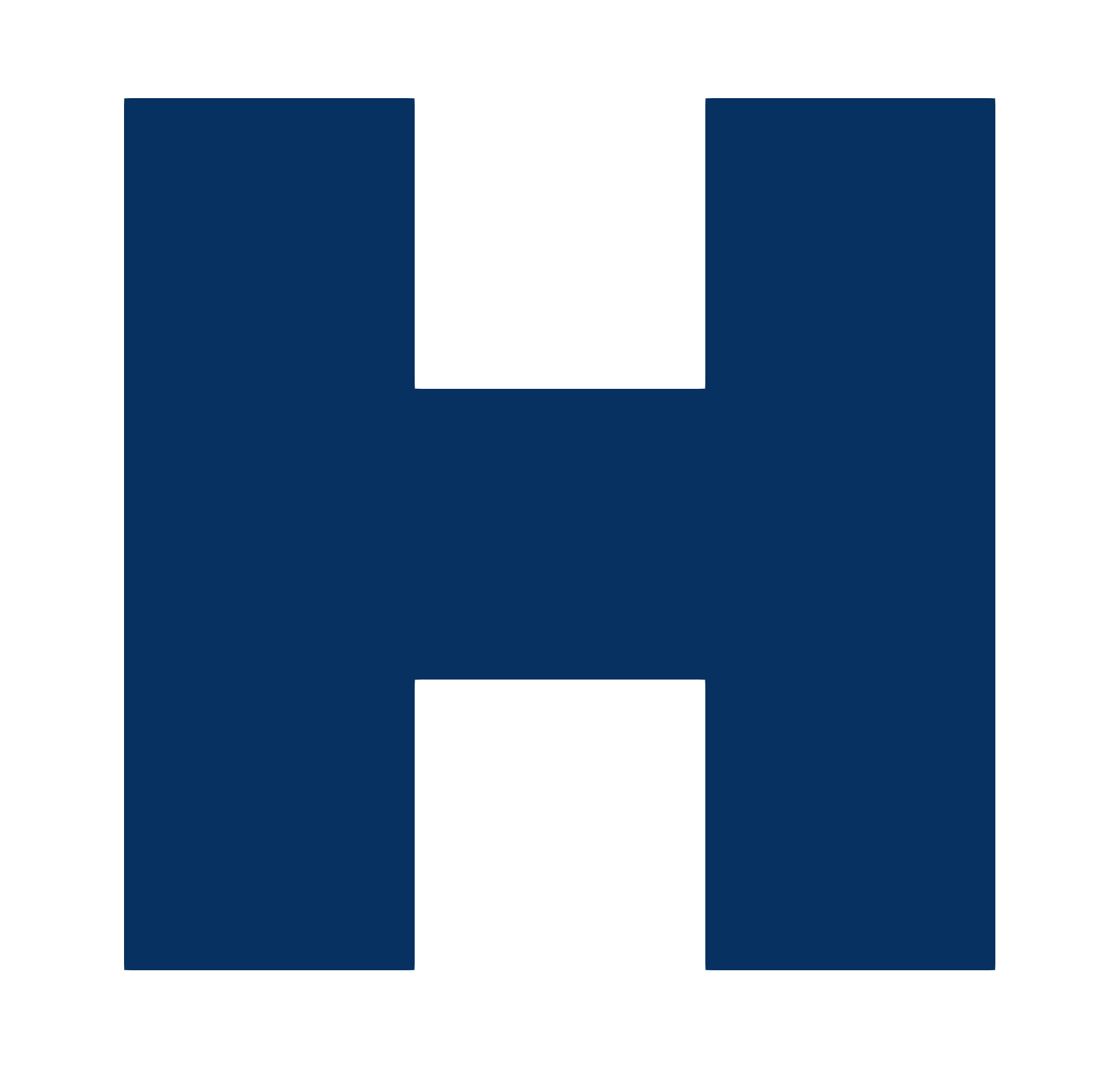}
  \hspace{.4cm}
  \includegraphics[height=6cm]{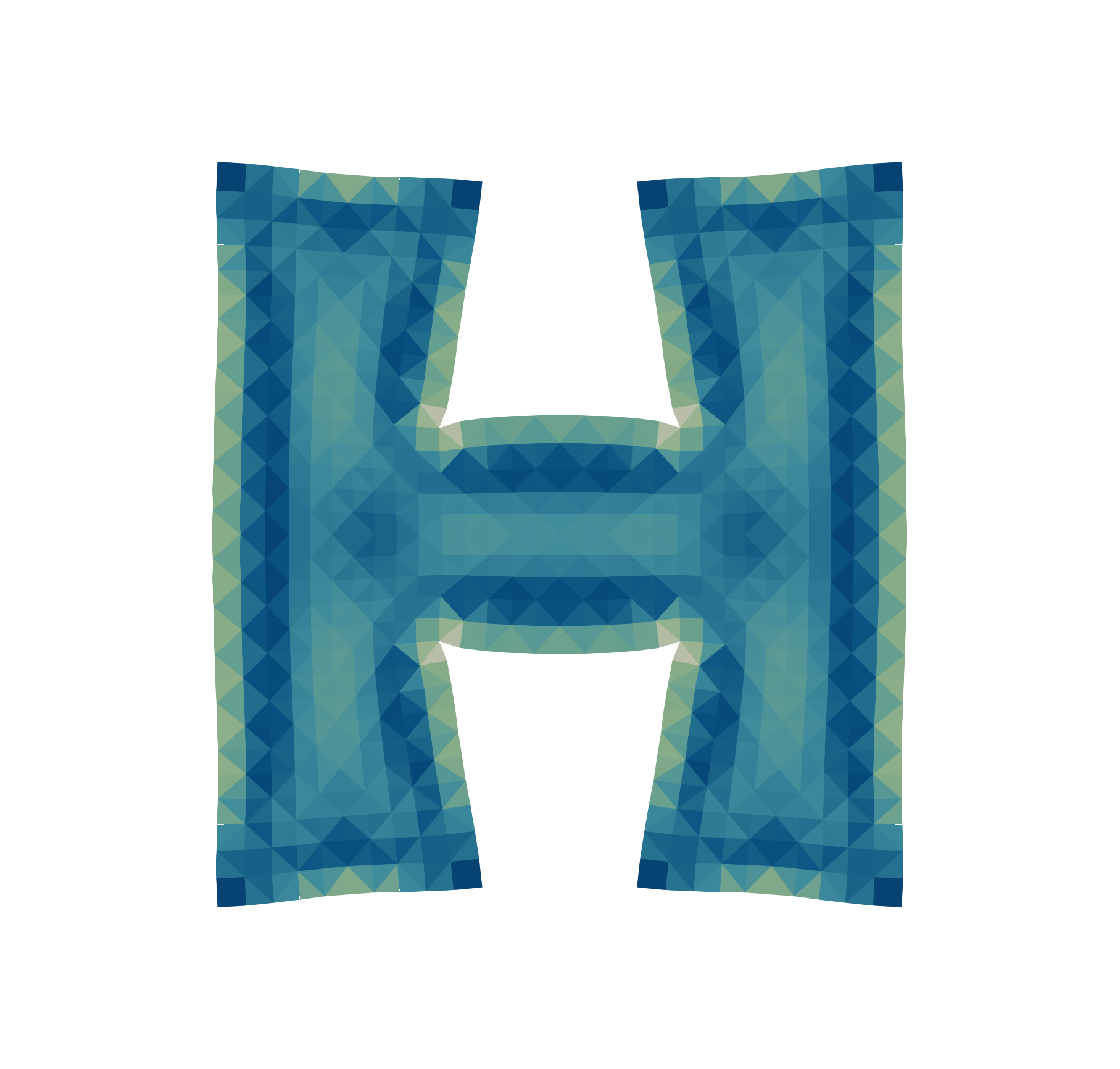}
  \hspace{.4cm}
  \includegraphics[height=6cm]{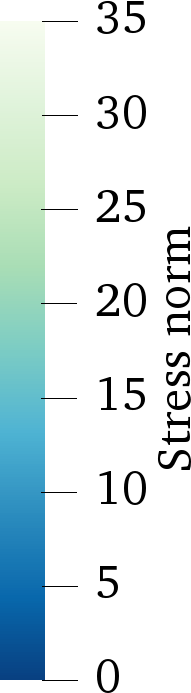}
  \caption{Spectral norm of initial Cauchy stress tensor in initial configuration (left)
  and residual (equilibrium) Cauchy stress tensor in equilibrium configuration (right)
  for thermoelastic body in forward simulation.}
  \label{fig:forward.snorm}
\end{figure}

In addition, Figure~\ref{fig:forward.snorm} shows the residual stresses in
the material, when starting from a stress-free state.

\begin{figure}[htbp]
  \centering
  \includegraphics[height=6cm]{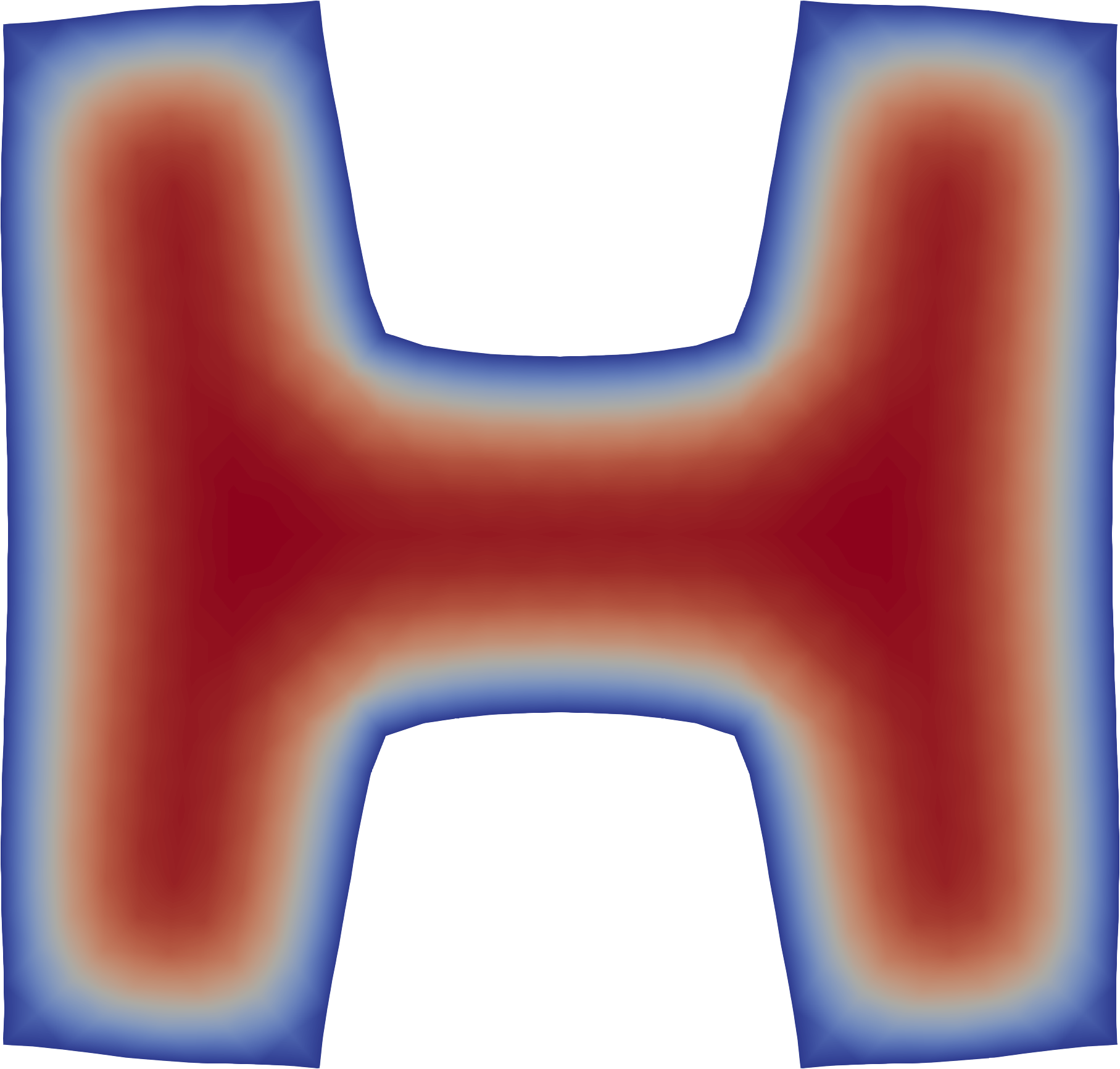}
  \hspace{.4cm}
  \includegraphics[height=6cm]{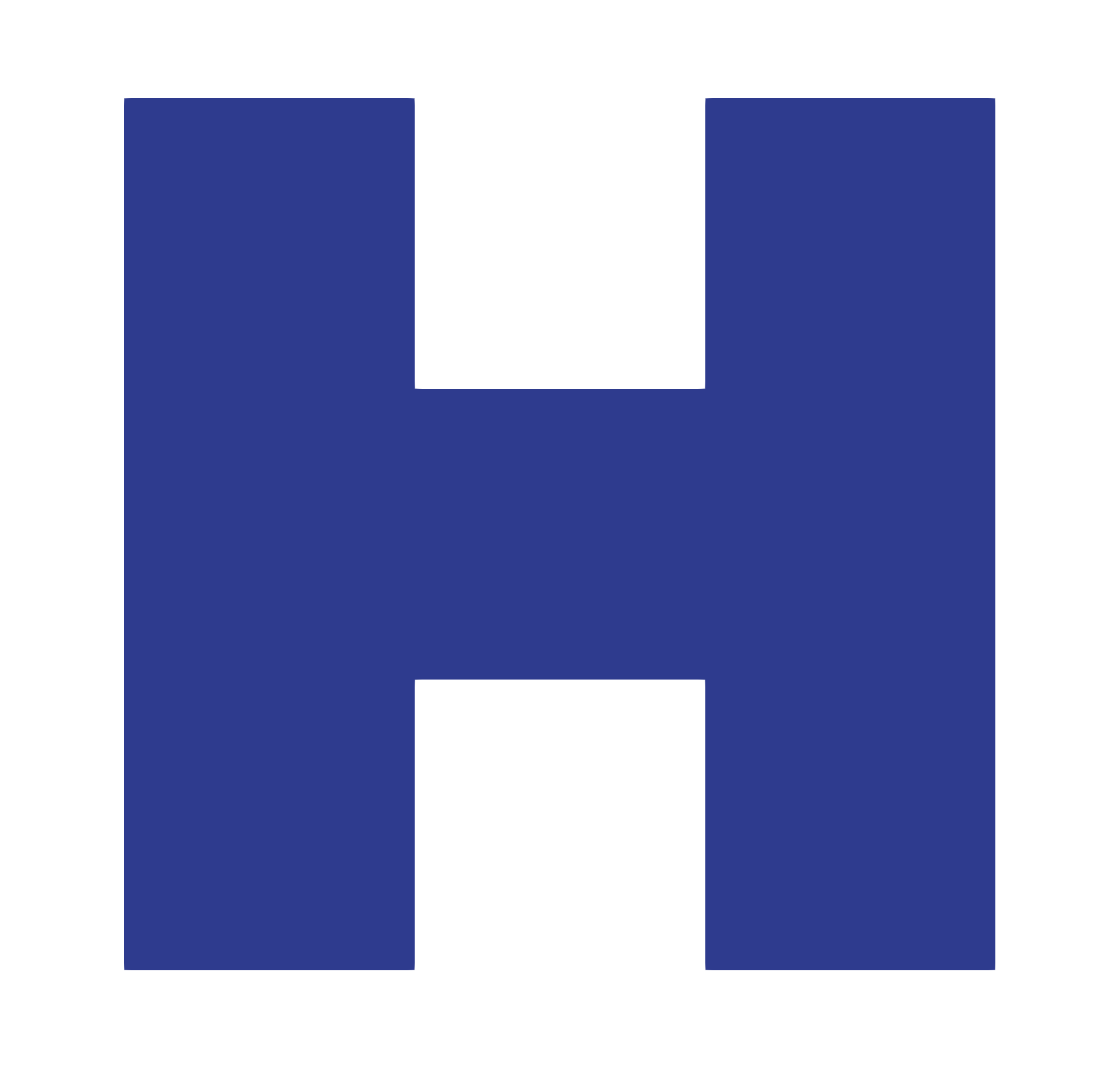}
  \hspace{.4cm}
  \includegraphics[height=6cm]{fig-thermal-legend-theta.png}
  \caption{Initial temperature in initial configuration (left)
  and equilibrium temperature in equilibrium configuration (right)
  for thermoelastic body in inverse simulation.}
  \label{fig:inverse.theta}
\end{figure}

The same inhomogeneous temperature field can be used as basis for the inverse
simulation. The left-hand side of Figure~\ref{fig:inverse.theta} shows the optimal
shape of the body at ejection time, under the assumption of minimal changes in
the temperature field, along with the adjusted temperature field.

\begin{figure}[htbp]
  \centering
  \includegraphics[height=6cm]{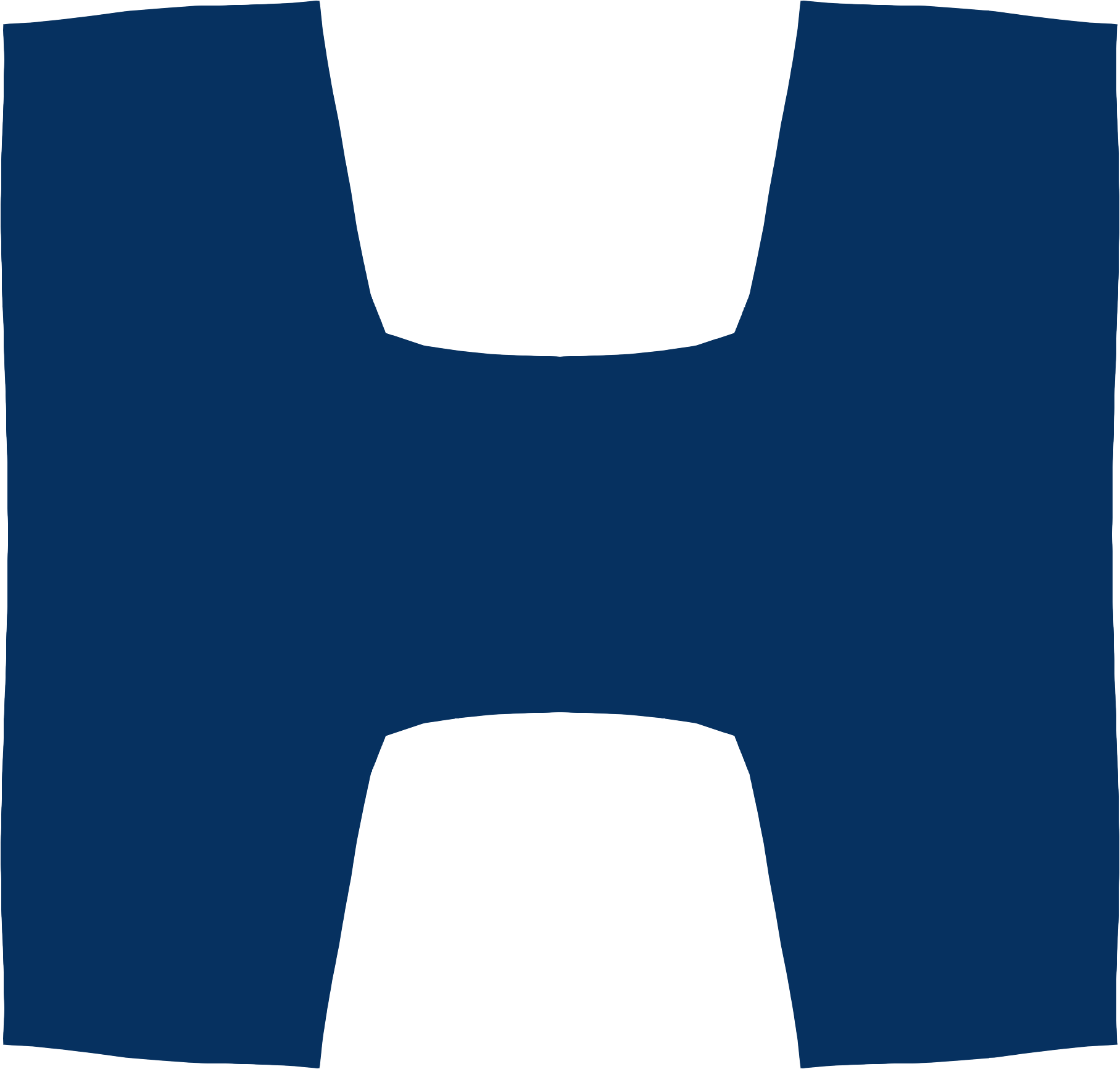}
  \hspace{.4cm}
  \includegraphics[height=6cm]{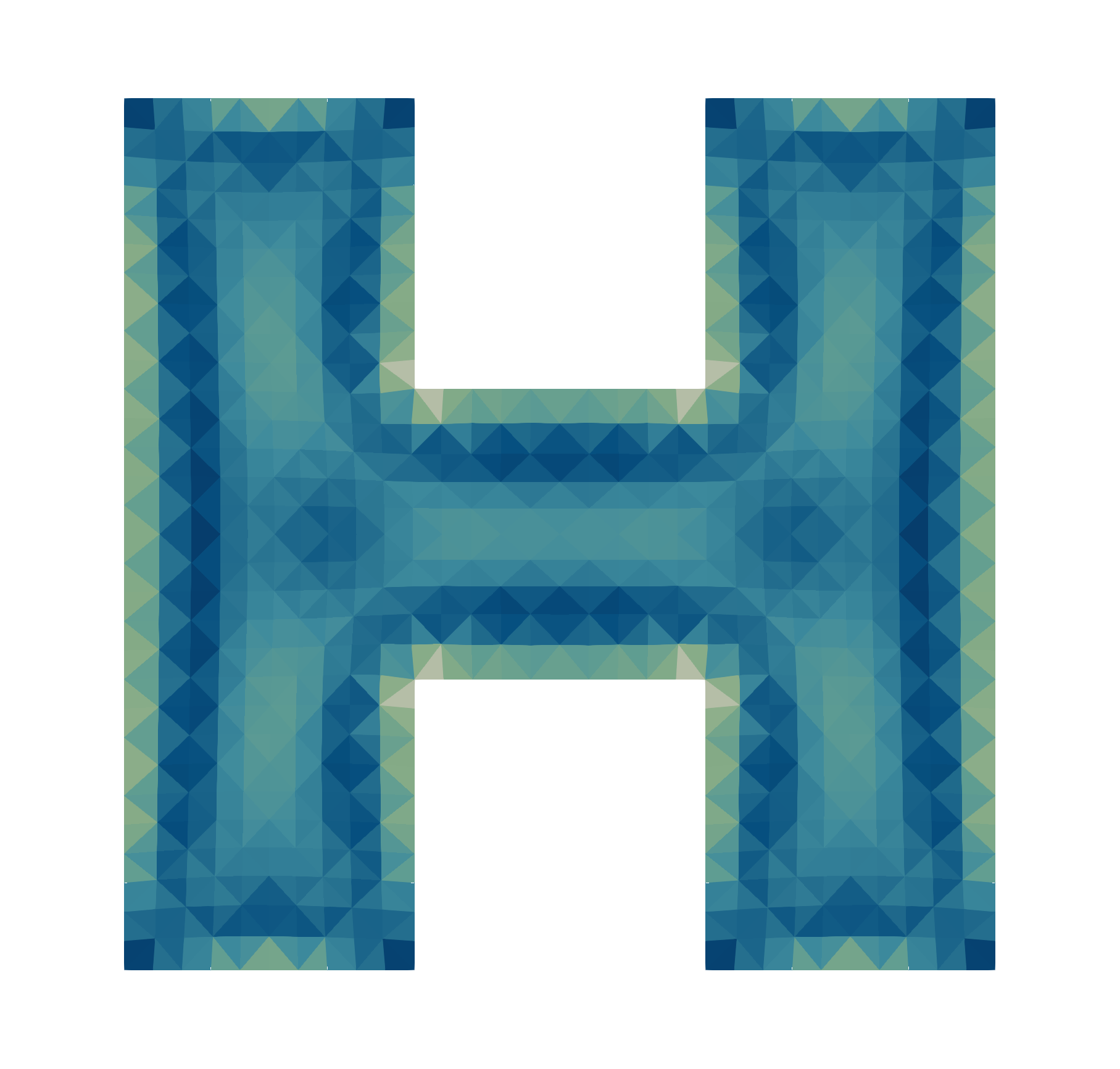}
  \hspace{.4cm}
  \includegraphics[height=6cm]{fig-thermal-legend-snorm.png}
  \caption{Spectral norm of initial Cauchy stress tensor in initial configuration (left)
  and residual (equilibrium) Cauchy stress tensor in equilibrium configuration (right)
  for thermoelastic body in inverse simulation.}
  \label{fig:inverse.snorm}
\end{figure}

Figure~\ref{fig:inverse.snorm} highlights how our simulation, while ensuring a
perfectly shaped result, does not prevent residual stresses in the body in any
way.

\section{Summary and Conclusion}

In this paper, we have shown a new method for the solution of inverse design
problems in thermoelasticity, including the formulation of the full, discretized
system of equations. Using numerical examples, we have also demonstrated the
usability of the method. The examples show that the method works and performs
according to intuitive expectations.

One of the key difficulties we dealt with was the prescription of temperature
and stress fields on an unknown shape. The proposed method of fitting these
fields, when given for a specific shape, into adjusted shapes in a way that
causes minimal changes is also demonstrated in the examples and
follows the expectations.

The examples show how, in contrast to many shape optimization problems, the
objective is achieved exactly. I.e., the desired shape is exactly obtained by
the method. This is an optimal result, but one should not forget that it comes at
the price of residual stresses in the body. However, since the stress field is
calculated as a by-product of the inverse method, subsequent analyses can be
performed on the stresses, e.g., if some material tolerances are known.

By choosing an inverse formulation rather than a shape optimization method,
we have avoided many issues that come with these methods, such as the need to
find a suitable objective function, obtain its derivatives, and choose a
low-dimensional parameterization of the unknown shape.
Overall,
this method has the potential of yielding better results than shape
optimization methods in less time. One should note that this is a direct
consequence of choosing fully elastic material laws. Any other laws that would
require a transient simulation cannot be handled by this method.

At this point, we are confident that this method is ready to be applied to
realistic problems in various fields.
As we have mentioned before,
the method was originally motivated by the cavity shape determination problem in injection molding.
Its development marks some considerable progress in this
direction, but some work remains to be done. Specifically, efficient methods for
the determination of the temperature and stress fields for provided cavity
shapes are still under development. Challenges such as the precise modeling of
polymer crystallization have to be faced to obtain optimum results.

\appendix

\section{Derivations}

\subsection{Reformulation of Constitutive Law Based on St.\ Venant-Kirchhoff}
\label{sec:reformstvk}

We have defined a constitutive law for thermoelasticity based on the
St.\ Venant-Kirchhoff material. This is given as a relationship between the second
Piola-Kirchhoff tensor $\secondpktot$ and the Green-St.\ Venant strain
tensor $\bm{E}$:

\begin{align}
  \label{eq:stvenkirtot}
  \secondpktot \; &= \; \lambda (\mathrm{tr} \, \bm{E}) \bm{I} + 2 \mu \bm{E}
- \alpha (\theta - \theta_0) \bm{I} \, .
\end{align}

Our differential equations require constitutive laws to be formulated for the
first Piola-Kirchhoff stress tensor $\firstpkref$ with respect to the reference
configuration.
The Green-St.\ Venant strain tensor $\bm{E}$ is defined as follows:

\begin{equation}
  \begin{aligned}
  \bm{E} \; &= \; \frac12 (\bm{C} - \bm{I})
  \; = \; \frac12 (\deftentot^T \deftentot - \bm{I})
  \; = \; \frac12 (\mathring{\bm{F}}^T \tilde{\bm{F}}^{-T} \bar{\bm{F}}^T
                   \bar{\bm{F}} \tilde{\bm{F}}^{-1} \mathring{\bm{F}} - \bm{I})
  \, .
  \end{aligned}
\end{equation}

The deformation gradient tensors $\tilde{\bm{F}}$ and $\bar{\bm{F}}$ can be
calculated directly from the fields $\vec{r}$ and $\vec{u}$, respectively
(see Figure~\ref{fig:configurations}).
However, since we assume the initial deformation to be provided in the form of
the left Cauchy-Green deformation tensor $\mathring{\bm{B}}$, the tensor
$\mathring{\bm{F}}$ is not available, and therefore we cannot calculate
$\bm{E}$.
In order to make Equation~(\ref{eq:stvenkirtot}) useful to us, we need to
reformulate it in such a way that the stress tensor $\firstpkref$ can
be calculated from the deformation tensors $\tilde{\bm{F}}$, $\bar{\bm{F}}$
and $\mathring{\bm{B}}$.
At this point it should be noted that the determinant $\mathring{J} :=
\mathrm{det} \, \mathring{\bm{F}}$ can be calculated from $\mathring{\bm{B}}$ as
$\mathring{J} = \sqrt{\mathrm{det} \, \mathring{\bm{B}}}$.

We will make use of the following identities \cite{Surana2015}:

\begin{alignat}{3}
  \bar{J} \boldsymbol{\sigma} \; &= \; \firstpkref \bar{\bm{F}}^T \, , \quad
  &J \boldsymbol{\sigma} \; &= \; \firstpktot \deftentot^T \, , \\
  \firstpkref \; &= \; \bar{\bm{F}} \secondpkref \, , \quad
  &\firstpktot \; &= \; \deftentot \secondpktot \, .
\end{alignat}

The first step is to find a relationship between the tensors $\firstpkref$ and
$\secondpktot$:

\begin{equation}
  \begin{aligned}
    \label{eq:relprefstot}
    && \firstpkref \bar{\bm{F}}^T \; &= \; \bar{J} \boldsymbol{\sigma} \\
    \Leftrightarrow \quad&&
    J \firstpkref \bar{\bm{F}}^T \; &= \; \bar{J} J \boldsymbol{\sigma} \\
    \Leftrightarrow \quad&&
    J \firstpkref \bar{\bm{F}}^T \; &= \; \bar{J} \firstpktot \deftentot^T \\
    \Leftrightarrow \quad&&
    \bar{J} \tilde{J}^{-1} \mathring{J} \firstpkref \bar{\bm{F}}^T
      \; &= \; \bar{J} \firstpktot \left(\bar{\bm{F}} \tilde{\bm{F}}^{-1}
               \mathring{\bm{F}}\right)^T \\
    \Leftrightarrow \quad&&
    \bar{J} \tilde{J}^{-1} \mathring{J} \firstpkref \bar{\bm{F}}^T
      \; &= \; \bar{J} \firstpktot
               \mathring{\bm{F}}^T \tilde{\bm{F}}^{-T} \bar{\bm{F}}^T \\
    \Leftrightarrow \quad&&
    \tilde{J}^{-1} \mathring{J} \firstpkref \; &= \; \firstpktot
               \mathring{\bm{F}}^T \tilde{\bm{F}}^{-T} \\
    \Leftrightarrow \quad&&
    \firstpkref \; &= \; \tilde{J} \mathring{J}^{-1} \deftentot \secondpktot
               \mathring{\bm{F}}^T \tilde{\bm{F}}^{-T} \, . \\
  \end{aligned}
\end{equation}

The following identities, which can be derived easily, will be helpful in the reformulation:

\begin{align}
  \deftentot \bm{I} \mathring{\bm{F}}^T \tilde{\bm{F}}^{-T}
    \; &= \; \bar{\bm{F}} \tilde{\bm{F}}^{-1} \mathring{\bm{B}}
             \tilde{\bm{F}}^{-T} \, , \\
  \deftentot \bm{C} \mathring{\bm{F}}^T \tilde{\bm{F}}^{-T}
    \; &= \; \bm{B} \bar{\bm{F}} \tilde{\bm{F}}^{-1} \mathring{\bm{B}} \tilde{\bm{F}}^{-T} \, , \\
  \mathrm{tr} \, \bm{E} \; &= \;
    \frac12 \mathrm{tr} \left( \bm{B} - \bm{I} \right)
    \, ,
\end{align}

where we have used $\bm{B} = \bar{\bm{F}} \tilde{\bm{F}}^{-1}
\mathring{\bm{B}} \tilde{\bm{F}}^{-T} \bar{\bm{F}}^T$.

Now, starting from Equation~(\ref{eq:relprefstot}) and using the other
relations, we can formulate the finished constitutive law:

\begin{equation}
  \begin{aligned}
    \firstpkref \; &= \;
      \tilde{J} \mathring{J}^{-1} \deftentot \secondpktot
      \mathring{\bm{F}}^T \tilde{\bm{F}}^{-T} \\
    &= \;
      \tilde{J} \mathring{J}^{-1} \deftentot
      \left( \lambda (\mathrm{tr} \, \bm{E}) \bm{I} + 2 \mu \bm{E}
    - \alpha (\theta - \theta_0) \bm{I} \right)
      \mathring{\bm{F}}^T \tilde{\bm{F}}^{-T} \\
    &= \;
      \tilde{J} \mathring{J}^{-1} \deftentot
      \left( \lambda (\mathrm{tr} \, \bm{E}) \bm{I} + \mu (\bm{C} - \bm{I})
    - \alpha (\theta - \theta_0) \bm{I} \right)
      \mathring{\bm{F}}^T \tilde{\bm{F}}^{-T} \\
    &= \;
      \tilde{J} \mathring{J}^{-1}
      \left( \frac{\lambda}2 \mathrm{tr} (\bm{B} - \bm{I}) \bm{I} + \mu (\bm{B} - \bm{I})
    - \alpha (\theta - \theta_0) \bm{I} \right)
      \bar{\bm{F}} \tilde{\bm{F}}^{-1} \mathring{\bm{B}}
             \tilde{\bm{F}}^{-T} \\
  \end{aligned}
\end{equation}

\subsection{Reformulation of Constitutive Law Based on Neo-Hooke}
\label{sec:reformneohooke}

The second constitutive law that we have mentioned is based on the Neo-Hooke
material. The Cauchy stress according to this law can be related to the left
Cauchy-Green deformation tensor $\bm{B}$ by

\begin{align}
  J \boldsymbol{\sigma} \; &= \; 2 D_1 J (J-1) \bm{I} + 2 C_1 J^{-\frac23} (\text{dev} \, \bm{B})
                   - \alpha (\theta - \theta_0) \bm{B} \, .
\end{align}

We need to formulate $\firstpkref$ in terms of $J \boldsymbol{\sigma}$:

\begin{equation}
  \begin{aligned}
    \label{eq:relprefcauchy}
    && \firstpkref \bar{\bm{F}}^T \; &= \; \bar{J} \boldsymbol{\sigma} \\
    \Leftrightarrow \quad&&
    J \firstpkref \bar{\bm{F}}^T \; &= \; \bar{J} J \boldsymbol{\sigma} \\
    \Leftrightarrow \quad&&
    \bar{J} \tilde{J}^{-1} \mathring{J} \firstpkref \bar{\bm{F}}^T
      \; &= \; \bar{J} J \boldsymbol{\sigma} \\
    \Leftrightarrow \quad&&
    \firstpkref
      \; &= \; \tilde{J} \mathring{J}^{-1} (J \boldsymbol{\sigma}) \bar{\bm{F}}^{-T}
    \, . \\
  \end{aligned}
\end{equation}

The law can now be inserted into the equation:

\begin{equation}
  \begin{aligned}
    \firstpkref
      \; &= \; \tilde{J} \mathring{J}^{-1} (J \boldsymbol{\sigma}) \bar{\bm{F}}^{-T} \\
    &= \; \tilde{J} \mathring{J}^{-1} \left(2 D_1 J (J-1) \bm{I} + 2 C_1 J^{-\frac23} (\text{dev} \, \bm{B})
                   - \alpha (\theta - \theta_0) \bm{B}\right) \bar{\bm{F}}^{-T} \, . \\
  \end{aligned}
\end{equation}

\section*{Declarations}

\subsection*{Competing interests}

The authors declare that they have no competing interests.

\subsection*{Funding}

The presented investigations were carried out at RWTH Aachen University within
the framework of the Collaborative Research Centre
SFB1120-236616214 "Bauteilpr\"azision durch Beherrschung von Schmelze und
Erstarrung in Produktionsprozessen" and funded by the Deutsche
Forschungsgemeinschaft e.V. (DFG, German Research Foundation). The sponsorship
and support is gratefully acknowledged.

\bibliography{references}

\begin{thebibliography}{10}
\expandafter\ifx\csname url\endcsname\relax
  \def\url#1{\texttt{#1}}\fi
\expandafter\ifx\csname urlprefix\endcsname\relax\def\urlprefix{URL }\fi
\expandafter\ifx\csname href\endcsname\relax
  \def\href#1#2{#2} \def\path#1{#1}\fi

\bibitem{Potsch2008}
G.~P{\"{o}}tsch, W.~Michaeli, {Injection Molding: An Introduction}, Carl Hanser
  Publishers, 2008.

\bibitem{Zwicke2017}
F.~Zwicke, M.~Behr, S.~Elgeti, {Predicting shrinkage and warpage in injection
  molding: Towards automatized mold design}, in: AIP Conference Proceedings,
  Vol. 1896, 2017.
\newblock \href {http://dx.doi.org/10.1063/1.5008119}
  {\path{doi:10.1063/1.5008119}}.

\bibitem{Liu1989}
D.~C. Liu, J.~Nocedal, {On the limited memory BFGS method for large scale
  optimization}, Mathematical Programming 45 (1989) 503--528.

\bibitem{Nobrega2008}
J.~M. N{\'{o}}brega, O.~S. Carneiro, A.~Gaspar-Cunha, N.~D. Goncalves, {Design
  of calibrators for profile extrusion - Optimizing multi-step systems},
  International Polymer Processing 23~(3) (2008) 331--338.
\newblock \href {http://dx.doi.org/10.3139/217.2148}
  {\path{doi:10.3139/217.2148}}.

\bibitem{Ettinger2004}
H.~J. Ettinger, J.~Sienz, J.~F.~T. Pittman, A.~Polynkin, {Parameterization and
  optimization strategies for the automated design of uPVC profile extrusion
  dies}, Structural and Multidisciplinary Optimization 28~(2-3) (2004)
  180--194.
\newblock \href {http://dx.doi.org/10.1007/s00158-004-0440-x}
  {\path{doi:10.1007/s00158-004-0440-x}}.

\bibitem{Elgeti2012}
S.~Elgeti, M.~Probst, C.~Windeck, M.~Behr, W.~Michaeli, C.~Hopmann, {Numerical
  shape optimization as an approach to extrusion die design}, Finite Elements
  in Analysis and Design 61 (2012) 35--43.
\newblock \href {http://dx.doi.org/10.1016/j.finel.2012.06.008}
  {\path{doi:10.1016/j.finel.2012.06.008}}.

\bibitem{Siegbert2015}
R.~Siegbert, N.~Yesildag, M.~Frings, F.~Schmidt, S.~Elgeti, H.~Sauerland,
  M.~Behr, C.~Windeck, C.~Hopmann, Y.~Queudeville, {Individualized production
  in die-based manufacturing processes using numerical optimization}, The
  International Journal of Advanced Manufacturing Technology 80~(5-8) (2015)
  851--858.

\bibitem{Schield1967}
R.~T. Schield, {Inverse deformation results in finite elasticity}, Zeitschrift
  f{\"{u}}r angewandte Mathematik und Physik ZAMP 18~(4) (1967) 490--500.
\newblock \href {http://dx.doi.org/10.1007/BF01601719}
  {\path{doi:10.1007/BF01601719}}.

\bibitem{Carlson1969}
D.~E. Carlson, T.~Shield, {Inverse deformation results for elastic materials},
  Zeitschrift f{\"{u}}r angewandte Mathematik und Physik ZAMP 20~(2) (1969)
  261--263.
\newblock \href {http://dx.doi.org/10.1007/BF01595564}
  {\path{doi:10.1007/BF01595564}}.

\bibitem{Govindjee1996}
S.~Govindjee, P.~A. Mihalic, {Computational methods for inverse finite
  elastostatics}, Computer Methods in Applied Mechanics and Engineering
  136~(96) (1996) 47--57.
\newblock \href {http://dx.doi.org/10.1016/0045-7825(96)01045-6}
  {\path{doi:10.1016/0045-7825(96)01045-6}}.

\bibitem{Govindjee1998}
S.~Govindjee, P.~A. Mihalic, {Computational methods for inverse deformations in
  quasi-incompressible finite elasticity}, International Journal for Numerical
  Methods in Engineering 43~(5) (1998) 821--838.

\bibitem{Hong2016}
W.~Hong, {Inverse Lagrangian formulation for the deformation of hyperelastic
  solids}, Extreme Mechanics Letters 9 (2016) 30--39.
\newblock \href {http://dx.doi.org/10.1016/j.eml.2016.04.009}
  {\path{doi:10.1016/j.eml.2016.04.009}}.

\bibitem{Yamada1998}
T.~Yamada, {Finite element procedure of initial shape determination for
  hyperelasticity}, Structural Engineering and Mechanics 6~(2) (1998) 173--183.
\newblock \href {http://dx.doi.org/10.12989/sem.1998.6.2.173}
  {\path{doi:10.12989/sem.1998.6.2.173}}.

\bibitem{Bazilevs2012}
Y.~Bazilevs, M.~C. Hsu, J.~Kiendl, D.~J. Benson, {A computational procedure for
  prebending of wind turbine blades}, International Journal for Numerical
  Methods in Engineering~(August 2011) (2013) 323--336.
\newblock \href {http://arxiv.org/abs/1010.1724} {\path{arXiv:1010.1724}},
  \href {http://dx.doi.org/10.1002/nme} {\path{doi:10.1002/nme}}.

\bibitem{Campbell2011}
R.~L. Campbell, {Fluid‐structure interaction and inverse design simulations
  for highly flexible turbomachinery.}, The Journal of the Acoustical Society
  of America 129~(4) (2011) 2385--2385.
\newblock \href {http://dx.doi.org/10.1121/1.3587740}
  {\path{doi:10.1121/1.3587740}}.

\bibitem{Fachinotti2008}
V.~D. Fachinotti, A.~Cardona, P.~Jetteur, {Finite element modelling of inverse
  design problems in large deformations anisotropic hyperelasticity},
  International Journal for Numerical Methods in Engineering 74~(6) (2008)
  894--910.
\newblock \href {http://arxiv.org/abs/1201.4903} {\path{arXiv:1201.4903}},
  \href {http://dx.doi.org/10.1002/nme.2193} {\path{doi:10.1002/nme.2193}}.

\bibitem{Albanesi2008}
A.~E. Albanesi, V.~D. Fachinotti, A.~Cardona, {Inverse Analysis of
  Large-Displacement Beams}, Mec{\'{a}}nica Computacional 27 (2008) 1049--1061.

\bibitem{Albanesi2013}
A.~E. Albanesi, M.~A. Pucheta, V.~D. Fachinotti, {A new method to design
  compliant mechanisms based on the inverse beam finite element model},
  Mechanism and Machine Theory 65 (2013) 14--28.
\newblock \href {http://dx.doi.org/10.1016/j.mechmachtheory.2013.02.009}
  {\path{doi:10.1016/j.mechmachtheory.2013.02.009}}.

\bibitem{Albanesi2011}
A.~E. Albanesi, {Inverse Design Methods for Compliant Mechanisms}, Ph.D. thesis
  (2011).

\bibitem{Albanesi2009}
A.~E. Albanesi, V.~D. Fachinotti, A.~Cardona, {Design of Compliant Mechanisms
  that Exactly Fit a Desired Shape}, Mec{\'{a}}nica Computacional 38~(38)
  (2009) 3191--3205.

\bibitem{Limache2011}
A.~C. Limache, {Inverse Shape Design of Deformable Structures and Deformable
  Wings}, Journal of Aircraft 48~(1) (2011) 157--165.
\newblock \href {http://dx.doi.org/10.2514/1.C001007}
  {\path{doi:10.2514/1.C001007}}.

\bibitem{Dennis2004}
B.~H. Dennis, G.~S. Dulikravich, S.~Yoshimura, {A Finite Element Formulation
  for the Determination of Unknown Boundary Conditions for Three-Dimensional
  Steady Thermoelastic Problems}, Journal of Heat Transfer 126~(1) (2004) 110.
\newblock \href {http://dx.doi.org/10.1115/1.1640360}
  {\path{doi:10.1115/1.1640360}}.

\bibitem{Dennis2011}
B.~Dennis, W.~Jin, {Application of the finite element method to inverse
  problems in solid mechanics}, International Journal of Structural Changes in
  Solids 3~(2) (2011) 11--21.

\bibitem{Zwicke2017a}
F.~Zwicke, S.~Eusterholz, S.~Elgeti, {Boundary-conforming free-surface flow
  computations: Interface tracking for linear, higher-order and isogeometric
  finite elements}, Computer Methods in Applied Mechanics and Engineering 326.
\newblock \href {http://dx.doi.org/10.1016/j.cma.2017.08.022}
  {\path{doi:10.1016/j.cma.2017.08.022}}.

\bibitem{Johnson1994}
A.~A. Johnson, T.~E. Tezduyar, {Mesh update strategies in parallel finite
  element computations of flow problems with moving boundaries and interfaces},
  Computer Methods in Applied Mechanics and Engineering 119~(1-2) (1994)
  73--94.
\newblock \href {http://dx.doi.org/10.1016/0045-7825(94)00077-8}
  {\path{doi:10.1016/0045-7825(94)00077-8}}.

\bibitem{Hughes2005}
T.~J.~R. Hughes, J.~A. Cottrell, Y.~Bazilevs, {Isogeometric analysis: CAD,
  finite elements, NURBS, exact geometry and mesh refinement}, Computer Methods
  in Applied Mechanics and Engineering 194~(39-41) (2005) 4135--4195.
\newblock \href {http://dx.doi.org/10.1016/j.cma.2004.10.008}
  {\path{doi:10.1016/j.cma.2004.10.008}}.

\bibitem{Surana2015}
K.~S. Surana, {Advanced Mechanics of Continua}, CRC Press, 2016.

\bibitem{Haupt2002}
P.~Haupt, J.~Wegner, {Continuum Mechanics and Theory of Materials}, Springer
  Science {\&} Business Media, 2013.
\newblock \href {http://dx.doi.org/10.1115/1.1451084}
  {\path{doi:10.1115/1.1451084}}.

\end{thebibliography}

\end{document}